\documentclass[11pt]{amsart}

\usepackage{amsmath,amssymb}
\usepackage{indentfirst}
\usepackage[all]{xy}
\newtheorem{tm}{Theorem}[section]
\newtheorem{lemma}[tm]{Lemma}

\newtheorem{defi}[tm]{Definition}

\newtheorem{cor}[tm]{Corollary}

\newtheorem{remark}[tm]{Remark}

\begin{document}

\newcommand{\R}{\mathbb{R}}
\newcommand{\N}{\mathbb{N}}
\newcommand{\Z}{\mathbb{Z}}
\newcommand{\Q}{\mathbb{Q}}
\newcommand{\scf}{\mathcal{F}}
\newcommand{\scc}{\mathcal{C}}
\newcommand{\scs}{\mathcal{S}}
\newcommand{\scu}{\mathcal{U}}
\newcommand{\scv}{\mathcal{V}}
\newcommand{\col}{\!:\!}
\newcommand{\ra}{\rightarrow}
\newcommand{\ov}{\overline}
\newcommand{\bs}{\backslash}
\renewcommand{\mid}{\big\vert}
\newcommand{\al}{\alpha}
\newcommand{\da}{\delta}
\newcommand{\Da}{\Delta}
\newcommand{\e}{\varepsilon}
\newcommand{\et}{\emptyset}
\newcommand{\ga}{\gamma}
\newcommand{\vp}{\varphi}
\newcommand{\Sa}{\Sigma}
\newcommand{\sa}{\sigma}
\newcommand{\la}{\lambda}
\renewcommand{\mod}{\mathrm{mod}}

\newcommand{\mesh}{\mathop{\mathrm{mesh}}}
\newcommand{\ord}{\mathop{\mathrm{ord}}}
\newcommand{\card}{\mathop{\mathrm{card}}}
\newcommand{\diam}{\mathop{\mathrm{diam}}}
\newcommand{\im}{\mathop{\mathrm{im}}}
\newcommand{\St}{\mathop{\mathrm{St}}}
\newcommand{\ANR}{\mathop{\mathrm{ANR}}}
\newcommand{\EW}{\mathop{\mathrm{EW}}}
\renewcommand{\int}{\mathop{\mathrm{int}}}
\newcommand{\inv}{^{-1}}
\newcommand{\Hom}{\mathop{\mathrm{Hom}}}
\newcommand{\Sh}{\mathop{\mathrm{Sh}}}
\newcommand{\UV}{\mathop{\mathrm{UV}}}
\newcommand {\Tor}{\mathop{\mathrm{Tor}}}


\title{Simultaneous $\Z / p$-acyclic resolutions of expanding sequences}

\author{Leonard R.\ Rubin}

\author{Vera Toni\'{c}$^{\ \ast}$}
\thanks{$^{\ast}$Part of this paper was written while the second author was a post-doctoral fellow at Nipissing University, Department of Computer Science and Mathematics, North Bay, Ontario, Canada}

\address{Department of Mathematics\\
University of Oklahoma\\
601 Elm Ave, room 423\\
Norman, Oklahoma 73019\\
USA}
\email{lrubin@ou.edu}

\address{Department of Mathematics\\
Ben Gurion University of the Negev\\
P.O.B.\ 653\\
Be'er Sheva 84105, Israel}
\email{vera.tonic@gmail.com}

\date{26 January, 2013}




\maketitle 


\markleft{ L.\ Rubin, V.\ Toni\'{c}}
\markright{Simultaneous $\Z / p$-acyclic resolutions}

\begin{abstract}
We prove the following theorem.

\noindent \textbf{Theorem}: \ \textit{Let
$X$ be a nonempty compact metrizable space, let $l_1\leq l_2\leq \ldots$ be a sequence in $\N$, and let $X_1\subset
X_2\subset\dots$ be a sequence of nonempty closed subspaces of $X$ such
that for each $k\in\N$, $\dim_{\Z/p} X_k\leq l_k$. Then there
exists a compact metrizable space $Z$, having closed subspaces
$Z_1\subset Z_2\subset\dots$, and a (surjective) cell-like map
$\pi:Z\ra X$, such that for each $k\in\N$,
\begin{enumerate}
\item[(a)]  $\dim Z_k\leq l_k$,
\item[(b)]  $\pi(Z_k)=X_k$, and
\item[(c)]  $\pi\vert_{Z_k}:Z_k\ra X_k$ is a $\Z/p$-acyclic map.
\end{enumerate}
Moreover, there is a sequence $A_1\subset A_2\subset\dots$ of
closed subspaces of $Z$ such that for each $k$, $\dim A_k\leq l_k$,
$\pi\vert_{A_k}:A_k\ra X$ is surjective, and for $k\in\N$,
$Z_k\subset A_k$ and $\pi\vert_{A_k}:A_k\ra X$ is a $\UV^{l_k-1}$-map.
}

\vspace{1mm}

It is not required that $X=\bigcup_{k=1}^\infty X_k$ or that $Z=\bigcup_{k=1}^\infty Z_k$.
This result generalizes the $\Z/p$-resolution theorem of A.\ Dranishnikov and runs parallel to a similar theorem
of S.~Ageev, R.~Jim\'{e}nez, and the first author, who studied the situation where the group was $\Z$.
\end{abstract}


\section{Introduction}

The goal of this paper is to prove the following theorem.

\begin{tm} \label{T1.1}
Let $X$ be a nonempty compact metrizable space, let $l_1\leq l_2\leq \ldots$ be a sequence in $\N$, and let $X_1\subset X_2\subset\dots$ be a sequence of nonempty closed subspaces of $X$ such
that for each $k\in\N$, $\dim_{\Z/p} X_k\leq l_k$. Then there
exists a compact metrizable space $Z$, having closed subspaces
$Z_1\subset Z_2\subset\dots$, and a (surjective) cell-like map
$\pi:Z\ra X$, such that for each $k\in\N$,
\begin{enumerate}
\item[(a)]  $\dim Z_k\leq l_k$,
\item[(b)]  $\pi(Z_k)=X_k$, and
\item[(c)]  $\pi\vert_{Z_k}:Z_k\ra X_k$ is a $\Z/p$-acyclic map.
\end{enumerate}
Moreover, there is a sequence $A_1\subset A_2\subset\dots$ of
closed subspaces of $Z$ such that for each $k$, $\dim A_k\leq l_k$,
$\pi\vert_{A_k}:A_k\ra X$ is surjective, and for $k\in\N$,
$Z_k\subset A_k$ and $\pi\vert_{A_k}:A_k\ra X$ is a $\UV^{l_k-1}$-map.
\end{tm}

The second Section will contain some technical results necessary for the proof of Theorem~\ref{T1.1}, and the proof will be described in the third Section. 

In Section 4 we will outline a proof of a case of Theorem~\ref{T1.1} that was suggested to us by an anonymous referee. Unfortunately, this technique cannot be used to prove the most difficult cases of Theorem~\ref{T1.1}, nor does it have the potential for generalization for those groups $G$ whose resolutions require a domain space of dimension $n+1$, if the range space had $\dim_G \leq n$ (\cite{Le}).
For example, the theorem that follows is an immediate consequence of Theorem~\ref{T1.1}, but it cannot be proven using the technique described in Section 4.

\begin{tm}
 Let $n\in \N$ and let $(X_i)$ be a sequence of (not necessarily nested) closed subsets
of the Hilbert cube $Q$ with $\dim_{\Z/p}X_i\leq n$ for all $i$.  Then
there exists a compact metrizable space $Z$, a cell-like
map $\pi:Z\to Q$, and a sequence $(Z_i)$
of closed subsets of $Z$ such that  $\forall i$,
\begin{enumerate}
\item[(a)] $\dim Z_i\leq n$, and 
\item[(b)] $\pi|_{Z_i}:Z_i\to X_i$ is a surjective
$\Z/p$-acyclic map.
\end{enumerate}
\end{tm}

This theorem provides a cell-like resolution
of the Hilbert cube $Q$ and simultaneously $\Z/p$-acyclic resolutions over any
F$_\sigma$-collection whatsoever of such $X_i$.

\vspace{2mm}

Let us proceed by explaining some terms that might be unfamiliar to the reader.
Basic facts about cell-like spaces and maps can be found in \cite{Da}. 
A \emph{map} $\pi: Z \ra X$ is called \emph{cell-like} if for each $x\in X$, $\pi^{-1}(x)$ has the shape of a point.  
To detect that a compact metrizable space $Y$ has the shape of a point, it is sufficient to prove that there is an inverse sequence $(Z_i,p_i^{i+1})$, of compact metrizable spaces $Z_i$, whose limit is homeomorphic to $Y$ and such that for infinitely many $i \in \N$, $p_i^{i+1}:Z_{i+1} \ra Z_i$ is null-homotopic. It is also sufficient to show that every map of $Y$ to a CW-complex is null-homotopic.

\vspace{1mm}

A map $\pi :Z \ra X$ between topological spaces is called $G$-\emph{acyclic} (\cite{Dr}) if
all its fibers $\pi^{-1}(x)$ have trivial reduced \v{C}ech
cohomology with respect to a given abelian group $G$, or, equivalently,  every map
$f: \pi^{-1}(x) \ra K(G,n)$ is null-homotopic.
Note that a map $\pi :Z \ra X$ being cell-like implies that $\pi$ is also $G$-acyclic.

To  detect that a compact metrizable space $Y$ has trivial reduced \v{C}ech
cohomology with respect to the group $G$, it is sufficient to prove that there is an inverse sequence $(Z_i,p_i^{i+1})$ of compact polyhedra $Z_i$ whose limit is homeomorphic to $Y$, such that for infinitely many $i\in \N$, the map $p_i^{i+1}: Z_{i+1} \ra Z_i$ induces the zero-homomorphism
of cohomology groups $\mathrm{H}^m (Z_i;G) \ra \mathrm{H}^m (Z_{i+1};G)$, for all $m \in \N$.

\vspace{1mm}

A map $\pi:Z\ra X$ is called a $\UV^k$-\textit{map} (\cite{Da})
if each of its fibers has property $\UV^k$.  This means that
each embedding $\pi^{-1}(x)\hookrightarrow A$ into an ANR $A$ has
property $\UV^k$: for every $0\leq r\leq k$ and every neighborhood $U$
of $\pi^{-1}(x)$ in $A$, there exists a neighborhood $V$ of $\pi^{-1}(x)$ in $U$
such that every map of $S^r$ into $V$ is null-homotopic in
$U$. In order to prove that $\pi$ is a $\UV^k$-map, it is sufficient to show that, for all $x\in X$,
there is an inverse sequence $(Z_i,p_i^{i+1})$ of compact polyhedra $Z_i$, whose limit is homeomorphic to $\pi^{-1}(x)$ and such that $\forall i \in \N$, if $0\leq r \leq k$, then any map $h:S^r \ra Z_i$ is null-homotopic.
It is well-known that cell-like
compacta have property $\UV^k$ for all $k$.

A map $g:X \ra \vert K\vert$ between a space $X$ and a polyhedron $|K|$ is called a $K$-\emph{modification} of a map $f:X\ra |K|$ if whenever $x\in X$ and $f(x)\in \sa$, for some $\sa \in K$, then $g(x)\in \sa$. This is equivalent to the following: whenever $x\in X$ and $f(x)\in \overset{\circ}\sa$, for some $\sa \in K$, then $g(x)\in \sa$.

\vspace{1mm}

The proof of Theorem~\ref{T1.1} uses some techniques developed by A.~Dranishnikov in the proof of the following theorem, which can be found as Theorem~8.7 in \cite{Dr}.

\begin{tm} \label{Dr} For every compact metrizable space $X$ with $\dim_{\Z/p} X$ $\leq n$,
there exists a compact metrizable space $Z$ and a surjective map
$\pi :Z \ra X$ such that $\pi$ is $\Z/p$-acyclic and $\dim Z \leq n$.
\end{tm}

We will show in Remark~\ref{R-Dr} that our
Theorem~\ref{T1.1} is a generalization of this theorem.
Dranishnikov used Edwards--Walsh complexes and resolutions, and so shall we.

The following definition of Edwards--Walsh complexes (EW-complexes) and resolutions, as well as results about them, can be found in \cite{Dr}, \cite{DW} or \cite{KY}. For $G=\Z$, these resolutions were formally formulated in \cite{Wa}.

\begin{defi}\label{EWDef}
Let $G$ be an abelian group,  $n\in\N$ and  $L$ a simplicial complex. An \emph{Edwards--Walsh resolution} of $L$ in dimension $n$ is a pair $(\EW(L,G,n),\omega)$ consisting of a \emph{CW}-complex $\EW(L,G,n)$ and a combinatorial map $\omega:\EW(L,G,n)\ra \vert L\vert$ (that is, $\omega^{-1} (\vert L'\vert)$ is a subcomplex, for each subcomplex $L'$ of $L$) such that:
\begin{enumerate}
\item[(i)] $\omega^{-1} (\vert L^{(n)}\vert)=\vert L^{(n)}\vert$ and $\omega\vert_{\vert L^{(n)}\vert}$ is the identity map of $\vert L^{(n)}\vert$ onto itself,
\item[(ii)] for every simplex $\sa$ of $L$ with $\dim \sa>n$, the preimage $\omega^{-1}(|\sa|)$ is an Eilenberg--MacLane space of type $(\bigoplus G,n)$, where the sum $\bigoplus G$ is finite, and
\item[(iii)] for every subcomplex $L'$ of $L$ and every map $f:\vert L'\vert \ra K(G,n)$, the composition $f\circ \omega\vert_{\omega^{-1} (\vert L'\vert)}:\omega^{-1} (\vert L'\vert)\ra K(G,n)$ extends to a map $F:\EW(L,G,n)\ra K(G,n)$.
\end{enumerate}
\end{defi}

We usually refer to the CW-complex $\EW(L,G,n)$ as an \emph{Edwards--Walsh complex}, and to the map $\omega$ itself as an \emph{Edwards--Walsh projection}.

\begin{remark}\label{R-EW}
 Let $L'$ be a subcomplex of $L$, $K$ be the subcomplex $\omega^{-1}(|L'|)$ of $EW(L,G,n)$, and
$\omega_{L'}=\omega|_{\omega^{-1}(|L'|)}:\omega^{-1}(|L'|)\to|L'|$. Then $(K,\omega_{L'})$ is an Edwards--Walsh resolution of the form $(EW(L',G,n),\omega_{L'})$.\\
\end{remark}
A discussion about the existence of Edwards--Walsh resolutions, as well as their construction, can be found in \cite{Dr}, \cite{DW}, \cite{KY}, \cite{Wa}.
For our needs, it is enough to know that when $G$ is either $\Z$ or $\Z/p$, Edwards--Walsh resolutions exist for any simplicial complex $L$. 

In particular, we shall briefly describe the construction of $(\EW(L,\Z/p,n),\omega)$ for a finite-dimensional simplicial complex $L$. If $\dim L\leq n$, define $\EW(L,\Z/p,n)=L^{(n)}=L$, and $\omega=id_L$.
If $\dim L= n+1$, we start with $L^{(n)}$ and the identity map $id_{L^{(n)}}$, and proceed by building a $K(\Z/p, n)$ on $\partial\sa$, for each $(n+1)$-simplex $\sa$ of $L$, and we build $\omega$ by extending $\partial\sa\hookrightarrow\sa$ over this newly attached $K(\Z/p,n)$. In this way, $\omega^{-1}(|\sa|)=K(\Z/p,n)$, $\forall$ $(n+1)$-simplex $\sa$ of $L$.

If $\dim L > n+1$, then we shall distinguish the cases $n\geq 2$ and $n=1$. In both of these cases our construction is inductive.

If $n\geq 2$ and $\dim L > n+1$, then the skeleton $L^{(n+1)}$ is dealt with as described above, i.e., by attaching a $K(\Z/p,n)$ to $\partial\sa$, for each $(n+1)$-simplex $\sa \in L^{(n+1)}$. This represents the basis of our inductive construction.
For the step of our inductive construction, let $k >n+1$. Then for any $k$-simplex $\sa$ of $L$, we have that $\pi_n(\omega^{-1}(|\partial\sa|))=\bigoplus\Z/p$, where this sum is finite. So $\omega^{-1}(|\sa|)$ will be obtained from $\omega^{-1}(|\partial\sa|)$ by attaching cells of $\dim \geq n+2$ in order to kill off the higher homotopy groups of $\omega^{-1}(|\partial\sa|)$, and achieve that $\omega^{-1}(|\sa|)=K(\bigoplus \Z/p,n)$.

If $n=1$ and $\dim L > 2$,  then the $2$-skeleton $L^{(2)}$ is dealt with as described above, that is, by attaching a $K(\Z/p,1)$ to $\partial\sa$, for each $2$-simplex $\sa \in L^{(n+1)}$. 
To be more precise, this means attaching a $2$-cell using a map of degree $p$ from the boundary of the $2$-cell to $\partial \sa$, for every $2$-simplex $\sa$ of $L$, and then proceeding by attaching cells of $\dim \geq 3$ to form a $K(\Z/p,1)$ on top of each of these Moore spaces. However, the above mentioned $2$-cells are not the only ones that get attached here; we will have to attach more of these.
Namely, when  $k> 2$, then for any $k$-simplex $\sa$ of $L$, there will be $2$-cells $\gamma\subset \omega^{-1}(|\sa|)\setminus \omega^{-1}(|\partial \sa|)$, attached by a map $\partial\gamma \ra \omega^{-1}(|\partial\sa|)$ representing a commutator in $\pi_1(\omega^{-1}(\partial\sa))$. This is to ensure that $\pi_1(\omega^{-1}(|\sa|))=\bigoplus \Z/p$.
We proceed by attaching cells of dimension $\geq 3$ to achieve that $\omega^{-1}(|\sa|)=K(\bigoplus \Z/p,1)$.

\vspace{2mm}

The following fact is proven in \cite{Dr} as Lemma~8.1, and (iv$_{\Z/p}$) is clear from our construction above.

\begin{lemma}\label{L8.1}
For the groups $\Z$ and $\Z/p$, for any $n\in\N$ and for any simplicial complex $L$, there is an Edwards--Walsh resolution  $\omega:\EW(L,G,n)\ra \vert L\vert$ with the additional property for $n>1$:
\begin{enumerate}
\item[(iv$_{\Z}$)]  the $(n+1)$-skeleton of $\EW (L,\Z,n)$ is equal to $L^{(n)}$;
\item[(iv$_{\Z/p}$)]  the $(n+1)$-skeleton of $\EW (L,\Z/p,n)$ is obtained from $L^{(n)}$ by attaching $(n+1)$-cells by a map of degree $p$ to the boundary $\partial \sa$, for every $(n+1)$-dimensional simplex $\sa$.
\end{enumerate}
\end{lemma}

Here are some other properties following from the construction of Edwards-Walsh complexes for the groups $\Z/p$.


\begin{remark}\label{No of summands}
Let $L$ be a simplicial complex, let $\sa$ be any simplex of $L$ with $\dim \sa>n$, and let $(\EW(L, \Z/p,n), \omega)$ be an Edwards-Walsh resolution of $L$.
According to Remark~\ref{R-EW}, $\omega^{-1}(\vert \sa \vert)=\EW (\sa, \Z/p, n)$ and from the construction of $\EW(L,\Z/p,n)$, we have that
  the number of summands in $\pi_n(\omega^{-1}(|\sa |))\cong \bigoplus \Z/p$ is less than or equal to the number of the $(n+1)$-faces of $\sa$. 
\end{remark}

From this Remark and our construction, we get:

\begin{cor}\label{simplex}
Let $\sa$ be a simplex with $\dim \sa>n$, taken as a simplicial complex, and let $(\EW(\sa, \Z/p,n), \omega)$ be an Edwards-Walsh resolution of $\sa$.
 Then
\begin{enumerate}
\item[(I)] $H_n(|\sa^{(n)}|)\cong \bigoplus_1^r \Z$, and
\item[(II)] $H_n(\EW(\sa,\Z/p,n))\cong \bigoplus_1^r \Z/p$, 
\end{enumerate}
where $r$ $\leq$ the number of all $(n+1)$-faces of $\sa$.
Moreover,
\begin{enumerate}
\item[(III)] we can choose $\tau_1, \ldots ,\tau_r$ to be some $(n+1)$-faces of $\sa$ so that the images $h_1 ,\ldots , h_r$ of the generators  of $H_n(\partial\tau_1), \ldots, H_n(\partial\tau_r)$, induced by the inclusions $\partial\tau_i \hookrightarrow  \sa^{(n)}$, form a basis of $H_n(|\sa^{(n)}|)$. Then if $q_1, \ldots , q_r$ are the images of the generators of $H_n(\partial\tau_1), \ldots, H_n(\partial\tau_r)$, induced by the inclusions $\partial\tau_i \hookrightarrow \EW(\sa,\Z/p,n)$, and $\la_\ast =H_n(\la)$ is induced by the inclusion $\la:\sa^{(n)}\hookrightarrow \EW(\sa,\Z/p,n)$, we get that $q_1=\la_\ast(h_1), \ldots , q_r=\la_\ast (h_r)$ form a basis of $H_n(\EW(\sa,\Z/p,n))$.
\end{enumerate}
\end{cor}

The following lemma is proven in \cite{Dr} as Lemma~8.2. It concerns (approximately) lifting maps through EW-complexes:
\begin{lemma}\label{L8.2}
Let $X$ be a compact
metrizable space with $\dim_G X \leq n$, and let $L$ be a finite simplicial complex. Then for every Edwards--Walsh resolution $\omega :
\EW(L,G,n) \to \vert L\vert$, and for every map $f:X \to \vert L\vert$, there
exists a map $f': X \to \EW(L,G,n)$ such that
\begin{enumerate}
\item[(i)] $f'\vert_{f^{-1}(\vert L^{(n)}\vert)}=f\vert_{f^{-1}(\vert L^{(n)}\vert)}$, and
\item[(ii)] $\omega\circ f'$ is an $L$-modification of $f$.
\end{enumerate}
\end{lemma}

\vspace{2mm}

Our primary construction will be done in the Hilbert cube $Q$ -- our space $X$ is compact metrizable, and $Q$ is universal for all compact metrizable spaces.

Let the Hilbert cube $\displaystyle Q=\prod_{i=1}^{\infty} I$ be
endowed with the metric $\rho$ such that if $x=(x_i)$, $y=(y_i)$,
then $\displaystyle \rho(x,y)=\sum_{i=1}^{\infty} \frac{\vert x_i-y_i\vert}{2^i}$. As usual, $I=[0,1]$.  For any $i\in\N$ it will be
convenient to write $Q=I^i\times Q_i$ in factored form.  In this
case, any subset $E$ of $I^i$ will always be treated as
$E\times\{0\}\subset Q$.  We shall use $p_i:Q\ra I^i$ for
coordinate projection.

\vspace{3mm}

In some of the proofs that follow we will use stability theory, about which more details can be found in \S VI.1 of \cite{HW}. Namely, we will use the consequences of Theorem VI.1. from \cite{HW}: if $X$ is a separable metrizable space with $\dim X \leq n$,  then for any map $f: X\ra I^{n+1}$ all values of $f$ are unstable. A point $y\in f(X)$ is called an  \emph{unstable value} of $f$ if for every $\da > 0$ there exists a map $g:X\ra I^{n+1}$ such that:
\begin{enumerate}
\item $d (f(x),g(x))<\da$ \ for every $x\in X$, and
\item$g(X)\subset I^{n+1}\setminus\{y\}$.
\end{enumerate}
Moreover, this map $g$ can be chosen so that $g=f$ on the complement of $f^{-1}(U)$, for any open neighborhood $U$ of $y$, and so that $g$ is homotopic to $f$ (see Corollary I.3.2.1 of \cite{MS}).

\vspace{2mm}

The following lemma is a form of the homotopy extension theorem
with control, and can be found in \cite{AJR} as Lemma~2.1.

\begin{lemma}\label{L2.1}
  Let $f:X\ra R$ be a map of a compact polyhedron $X$
to a space $R$, $X_0$ be a closed subpolyhedron of $X$, and
$\mathcal{U}$ be an open cover of $R$.  Suppose that $F:X_0\times
I\ra R$ is a $\mathcal{U}$-homotopy of $f\vert_{X_0}$. Then there
exists a $\mathcal{U}$-homotopy $H:X\times I\ra R$ of $f$ such
that $H\vert_{X_0\times I}=F:X_0\times I\ra R$.
\end{lemma}

\noindent\textbf{Notation.} We will use the following notation.  Let $x$ belong to a metric space
$X$ and let $\da>0$.  Then by $\overline N(x,\da)$ we shall mean the
closed $\da$-neighborhood of $x$ in $X$. Usually there will be no ambiguity, but notice that for $x\in Q$, $p_i(x)\in I^i$ 
so $\overline N(p_i(x),\da)$ will always refer to the
closed $\da$-neighborhood of $p_i(x)$ in $I^i$, even though $p_i(x)$ might also be contained in some subsets of $I^i$. If $\sa$ is a simplex in a triangulation $\tau$ of a polyhedron $P$, then $N(\sa,\da)$ will stand for the open $\da$-neighborhood of $\sa$ in $P$.

Whenever $(P_i,g_i^{i+1})$ is an
inverse sequence, $T_i\subset P_i$ and $g_i^{i+1}(T_{i+1})\subset T_i$ for
each $i$, then we shall write $(T_i,g_i^{i+1})$ for the induced inverse
sequence, using the same notation for the bonding maps as long as no confusion
can arise.

Whenever $P$ is a polyhedron, $\tau$ is a triangulation of $P$,
and $k\geq 0$, then $P^{(k)}$ will denote the subpolyhedron of $P$ triangulated by the $k$-skeleton of $\tau$, i.e., $P^{(k)}=|\tau^{(k)}|$.
If $R$ is a subpolyhedron of $P$ and we have to build an Edwards-Walsh complex  on $\tau|_R$, we will write $\EW(R,\Z/p,n)$ instead of $\EW(\tau|_R,\Z/p,n)$, to keep matters simpler.

\section{Technical lemmas}

The following type of result is a lemma which is technical, but
which will help us find certain maps and understand their fibers. This lemma can be found in
\cite{AJR} as Lemma~3.1.
Once the correct conditions are found on the construction of said maps, then
 Theorem~\ref{T1.1} will follow readily.

\begin{lemma}\label{L3.1}
 Suppose that for each $i\in\N$ we have selected
$n_i\in\N$, a compact subset $P_i\subset I^{n_i}$, $\da_i>0$,
$\e_i>0$, and a map $g_i^{i+1}:P_{i+1}\ra P_i$ so that:
\begin{enumerate}
 \item[(i)]  if $u$, $v\in Q$ and
$\rho(u,v)\leq\e_{i+1}$, then $\rho(p_{n_i}(u),p_{n_i}(v))<\da_i$,
\item[(ii)] $n_i<n_{i+1}$,
\item[(iii)] $\frac{9}{2^{n_i}}<\e_i$,
\item[(iv)] $\rho(g_i^{i+1}(x),p_{n_i}(x))<\da_i$ for all $x\in P_{i+1}$,
\item[(v)] $\da_i<\frac{1}{2^{n_i-1}}$, and
\item[(vi)] $P_{i+1}\times Q_{n_{i+1}}\subset P_i\times Q_{n_i}$.
\end{enumerate}

Put $\displaystyle X=\bigcap_{i=1}^{\infty} P_i\times Q_{n_i}$,
$\bold P=(P_i,g_i^{i+1})$, and $Z=\lim\mathbf P$.  Then for each
$\displaystyle z=(a_1,a_2,\dots)\in Z
\subset\prod_{i=1}^{\infty}P_i$, and associated sequence $(a_i)$
in $Q$, \begin{enumerate}
 \item[(a)]
$(a_i)$ is a Cauchy sequence in $Q$ whose limit lies in $X$, and
\item[(b)] the function $\pi:Z\ra X$ given by $ \displaystyle
\pi(z)= \lim_{i\ra\infty}(a_i)$ is continuous.
\end{enumerate}

Fix $x\in X$ and for each $i\in\N$, let $B_{x,i}=\overline
N(p_{n_i}(x),2\da_i)\cap P_i, B_{x,i}^\#=\overline
N(p_{n_i}(x),\e_i)\cap P_i$. Then,
\begin{enumerate} \item[(c)]
$B_{x,i}\subset B_{x,i}^\#$ and
 $g_i^{i+1}(B_{x,i+1}^\#)\subset B_{x,i}$.
\end{enumerate}

If we let $\bold P_x=(B_{x,i},g_i^{i+1})$ and $\bold P_x^\#=
(B_{x,i}^\#,g_i^{i+1})$, then, \begin{enumerate}
 \item[(d)]
$\lim\bold P_x=\lim\bold P_x^\#$, and \item[(e)]
$\pi^{-1}(x)=\lim\bold P_x$.
\end{enumerate}

In addition, suppose we are given, for each $i\in\N$, a closed
subspace $T_i\subset P_i$ in such a manner that
$g_i^{i+1}(T_{i+1})\subset T_i$. Put $\bold T=(T_i,g_i^{i+1})$ and
$Z'=\lim\bold T\subset Z$. For $x\in X$, let $S_{x,i}=B_{x,i}\cap
T_i$, $\bold T_x=(S_{x,i},g_i^{i+1})$; set $\tilde\pi=\pi\vert_{Z'}:Z'\ra X$. Then,
\begin{enumerate}
 \item[(f)]
$\tilde\pi^{-1}(x)=\lim\bold T_x$, and \item[(g)]  if
$S_{x,i}\neq\et$ for each $i$, then
$\tilde\pi^{-1}(x)\neq\emptyset$.
\end{enumerate}
\end{lemma}

A helpful diagram for Lemma~\ref{L3.1}:
\begin{displaymath}
\xymatrix{
\ldots &P_i  \ar[l] \ar@{_{(}->}[d]
&  P_{i+1}  \ar@{_{(}->}[d] \ar[l]^{p_{n_i}\vert}
\ar@/_/@{..>}[l]_{g_i^{i+1}} &\ \ldots \ar[l]&Z \ar@{-->}[d]^{\pi}\\
\ldots &P_i\times Q_{n_i}  \ar@{_{(}->}[l]&  P_i\times Q_{n_{i+1}}  \ar@{_{(}->}[l]&  \ \ldots \ar@{_{(}->}[l] & X
}
\end{displaymath}

Before proceeding, note that if $L$ is a simplicial
complex, $K$  a CW-complex, and $f : |L| \ra K$  a map, then we say that
$f$ is \emph{cellular} if it is cellular with respect to the CW-structure induced on
$|L|$ by $L$ and the given one of K, i.e., $f$ takes the (simplicial) $n$-skeleton of $L$ to the (CW) $n$-skeleton of $K$, $\forall n$.

The following Corollary is a  version of Corollary~3.2 from \cite{AJR}, adapted for the $\Z/p$-case.
When used (in the proof of the main theorem), $A_k$ can be replaced by $Z_k$ (not just by $A_k$ of Theorem~\ref{T1.1}).

\begin{cor}\label{C2.2} 
Suppose in Lemma~\ref{L3.1} that for each $i\in\N$, $P_i=\vert \tau_i\vert$ is a nonempty subpolyhedron of
$I^{n_i}$ having a triangulation $\widetilde\tau_i$, with a subdivision $\tau_i$ with $\mesh\tau_i<\da_i$, so that for every simplex $\ga$\ of $\widetilde \tau_i$, $\tau_i\vert_{\ga}$ is collapsible. 
Moreover, assume that  $g_i^{i+1}$ is a simplicial map (in particular, for all $k\geq 0$, $g_i^{i+1}(P_{i+1}^{(k)})\subset P_i^{(k)}$, where $\tau_{i+1}$ and $\tau_i$ are the relevant triangulations). 
Let $l_1\leq l_2\leq \ldots$ be a sequence in $\N$, and let
$$\bold T_k=(P_i^{(l_k)},g_i^{i+1}),
\,\,\mathrm{and}\,\,A_k=\lim\bold T_k.$$ Then
$A_1\subset
A_2\subset\dots$, and for each $k\geq 1$,
\begin{enumerate}
\item[(I)] $\dim A_k\leq l_k$ and $\pi\vert_{A_k}:A_k\ra X$ is surjective.
\end{enumerate}
Assume further that for each $x\in X$ and $i\in\N$, there is a contractible polyhedron $P_{x,i}$ which is the
closed star of a vertex in the triangulation $\widetilde\tau_i$, such that
$$B_{x,i}\subset P_{x,i}\subset B_{x,i}^\# .$$
 Then
\begin{enumerate}
 \item[(II)] $\pi:Z\ra X$ is a cell-like map, and
 \item[(III)]  for each $k\in\N$, $\pi\vert_{A_k}:A_k\ra X$ is a $\UV^{l_k-1}$-map.
\end{enumerate}
Suppose now that all of the above statements are true, and let
$k\in\N$. If for infinitely many indexes $i$ we have that
for all $x \in X$, $\omega \circ
\bar{f_i} (P_{x,i+1}) \subset P_{x,i}$, and $g_i^{i+1}\vert_{P_{x,i+1}} \simeq \omega \circ \bar{f_i}\vert_{P_{x,i+1}}$,
where $\omega :
\EW(P_i,\Z/p,l_k) \ra P_i$ is an Edwards--Walsh projection, and
$ \bar{f_i} : P_{i+1} \ra \EW(P_i,\Z/p,l_k)$ is a cellular
map, then
 \begin{enumerate}
\item[(IV)]  $\pi\vert_{A_k}:A_k\ra X$ is a $\Z/p$ -acyclic map.
\end{enumerate}
\end{cor}

Before showing the proof of Corollary~\ref{C2.2}, we will state and prove some lemmas which will be useful for its proof.

\begin{lemma}\label{algebra}
Let $n\in\N$, and let $P=\vert L\vert$  and $Q =\vert M \vert$ be compact polyhedra with $\dim P$, $\dim Q \geq n+1$. 
For any $(n+1)$-simplex $\tau_e$ of $M$, let $h_e$ and $q_e$ be the images of a generator of $H_n(\partial \tau_e)$ under the homomorphisms of $H_n(\partial \tau_e)$ induced by the inclusions
$\partial \tau_e \hookrightarrow \vert M^{(n)}\vert$ and $\partial\tau_e\hookrightarrow \EW(M,\Z/p,n)$, respectively.

Let $\mu$, $\nu$ and $\lambda$ be the inclusions as shown in the upcoming  diagram, and let $f:\vert L\vert \ra \EW (M, \Z/p, n)$ be a cellular map making this diagram commutative.

Moreover, let $M$ be such that:
\begin{enumerate}
\item[(I)] $H_n(| M^{(n)}|)\cong \bigoplus_1^r \Z$, and
\item[(II)] $H_n(\EW( M,\Z/p,n))\cong \bigoplus_1^r \Z/p$, 
\end{enumerate}
where $r$ $\leq$ the number of all $(n+1)$-simplexes of $ M$; and
\begin{enumerate}
\item[(III)] we can choose some $(n+1)$-simplexes $\tau_1, \ldots ,\tau_r$ of  $M$ so that $\{h_1 ,\ldots , h_r\}$ forms a basis of $H_n(| M^{(n)}|)$, and so that $\{q_1, \ldots , q_r\}$ forms a basis of $H_n(\EW(M,\Z/p,n))$.
\end{enumerate}

Then for any $(n+1)$-simplex $\sa \in L$, with $H_n(\partial \sa)=\langle g\rangle$, we have:
\begin{enumerate}
\item[(a)] $f\circ \nu\circ \mu$ is null-homotopic, so
\item[(b)] $H_n(f\vert_{\vert L^{(n)}\vert} \circ \mu)(g) = \sum_{e=1}^r \e_e h_e$,  
where $\e_e \equiv 0 \ (\mod \ p)$,  for $e\in\{1, \ldots , r\}$ .
\end{enumerate}
\end{lemma}

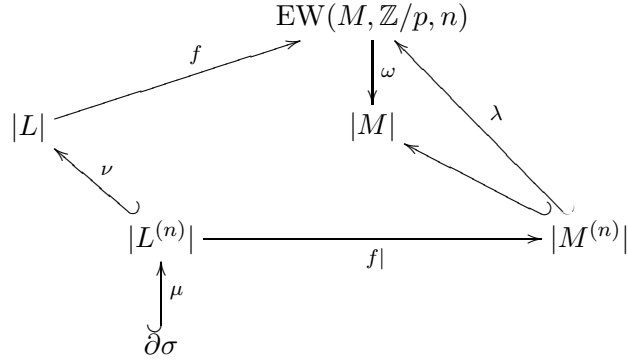
\begin{displaymath}
\xymatrix{
&&& \EW(M,\Z/p,n) \ar[d]^{\omega}   &\\
&\vert L\vert \ar[rru]^{f} && \vert M\vert &\\
&&\vert L^{(n)}\vert \ar@{_{(}->}[ul]_{\nu} \ar[rr]_{f\vert}  && \vert M^{(n)}\vert \ar@{_{(}->}[ul] \ar@{_{(}->}[uul]_{\la} \\
&& \partial \sa \ar@{^{(}->}[u]_{\mu}&&
}
\end{displaymath}

\noindent \emph{Proof}:
Since $\partial \sa$ is contained in $\sa$, which is contractible,
the inclusion $\nu \circ\mu :  \partial \sa \hookrightarrow \vert L\vert$ is null-homotopic.
Therefore $f\circ \nu\circ \mu$ is null-homotopic, so (a) is true.

To prove (b), notice that $f$ being a cellular map implies
$$f(\vert L^{(n)}\vert) \subset \EW(M,\Z/p,n)^{(n)} =
\vert M^{(n)}\vert.$$

It is clear that $f\circ \nu \circ \mu =
\la \circ f\vert_{\vert L^{(n)}\vert}\circ \mu$.
So (a) implies
$$0=H_{n}(f\circ \nu \circ \mu)(g)=H_{n}(\la \circ f\vert_{\vert L^{(n)}\vert}\circ
\mu)(g).$$

From (III) we get that 
$H_{n}(f\vert_{\vert L^{(n)}\vert}\circ
\mu)(g)=\sum_{e=1}^r \e_eh_e, \mathrm{\ for\ some\ } \e_e \in \Z$,  and therefore
$$H_{n}(\la \circ f\vert_{\vert L^{(n)}\vert}\circ
\mu)(g) = H_{n}(\la) (\sum_{e=1}^r \e_eh_e)=\sum_{e=1}^r \e_eq_e =0, $$
which means that $\e_e\equiv 0 \ (\mod \ p)$, $\forall e\in\{1, \ldots ,r\}$. $\square$

\vspace{2mm}

Some form of the following lemma was used by various authors.

\begin{lemma}\label{H_n(EW)}
Let $n\in\N$, $P=\vert \widetilde L\vert$ be a compact polyhedron with $\dim P \geq n+1$ and $\widetilde M$ be the closed star of a vertex from $\widetilde L^{(0)}$. Let $L$ be a subdivision of $\widetilde L$ such that for every simplex $\sa$ of $\widetilde L$, $L\vert_{\vert\sa\vert}$ is a collapsible simplicial complex. Let $M$ be the simplicial complex that $L$ induces on $|\widetilde M|$, i.e., $M=L\vert_{\vert \widetilde M\vert}$ (subdivided vertex star). 
 Then
\begin{enumerate}
\item[(I)] $H_n(| M^{(n)}|)\cong \bigoplus_1^r \Z$, and
\item[(II)] $H_n(\EW( M,\Z/p,n))\cong \bigoplus_1^r \Z/p$, 
\end{enumerate}
where $r$ $\leq$ the number of all $(n+1)$-simplexes of $ M$.
Moreover,
\begin{enumerate}
\item[(III)] we can choose $\tau_1, \ldots ,\tau_r$ to be some $(n+1)$-simplexes of $ M$ so that the images $h_1 ,\ldots , h_r$ of  the generators  of $H_n(\partial\tau_1), \ldots, H_n(\partial\tau_r)$, induced by the inclusions $\partial\tau_i \hookrightarrow   M^{(n)}$, form a basis of $H_n(| M^{(n)}|)$. Then if $q_1, \ldots , q_r$ are the images of the generators of $H_n(\partial\tau_1), \ldots, H_n(\partial\tau_r)$, induced by the inclusions $\partial\tau_i \hookrightarrow \EW( M,\Z/p,n)$, and $H_n(\la) $ is induced by the inclusion $\la: M^{(n)}\hookrightarrow \EW( M,\Z/p,n)$, we get that $q_1=H_n(\la)(h_1), \ldots , q_r=H_n(\la) (h_r)$ form a basis of $H_n(\EW(M,\Z/p,n))$.
\end{enumerate}
\end{lemma}
We will omit the proof to save space. On the way to proving this Lemma, one can first use Corollary~\ref{simplex} (containing the statement analogous to this one, but for a simplex) in order to prove analogous statements for a (non-subdivided) vertex star, and then for a subdivided simplex with a collapsible subdivision.

Then  Lemma~\ref{H_n(EW)} can be proven by first proving its statement for $\dim M =n+1$, and then, by induction, showing it is true for $\dim M =n+k+1$. The general step of induction would utilize another induction, on the number of $(n+k+1)$-simplexes of $\widetilde M$, as well as a Mayer-Vietoris sequence.  We used a collapsible subdivision on simplexes of $\widetilde M$ so that we could organize the process of induction. The information about the existence of subdivisions of a triangulation on a simplicial complex, in which a simplex with a new subdivision is still collapsible can be found in \cite{Gl}.

\begin{remark}\label{vertex star}
When $M$ is a subdivided vertex star from Lemma~\ref{H_n(EW)}, then Lemma~\ref{algebra} is true for $Q=|M|$ and $|M^{(n)}|$ is $(n-1)$-connected.
\end{remark}

\vspace{3mm}

\noindent \textsl{Proof of Corollary~\ref{C2.2}}: 
Surely $\dim A_k\leq l_k$. Let $x\in X.$
Apply Lemma~\ref{L3.1} with $T_i=P_i^{(l_k)}$ and $$S_{x,i}=B_{x,i}\cap P_i^{(l_k)}.$$ Then $\bold T$
becomes $\bold T_k$ and
$$Z'= \lim \bold T_k= \lim (P_i^{(l_k)},g_i^{i+1})=A_k.$$
Note that the representation of $X$ implies that $p_{n_i}(X)\subset P_i$, $\forall i\in\N$.
This fact, together with $\mesh\tau_i<\da_i$, can be used  to check that $B_{x,i}$ must contain a vertex of $\tau_i$, so $S_{x,i}\neq\et$.  Therefore  (g) of Lemma~\ref{L3.1} shows that (I) is true.

Part (c) of Lemma~\ref{L3.1} and the fact that
$B_{x,i}\subset P_{x,i}\subset B_{x,i}^{\#}$  $\forall i\in\N$,
show that $\forall i \in\N$, $g_i^{i+1}(P_{x,i+1})\subset P_{x,i}$, so
$\bold P_x':=(P_{x,i},g_i^{i+1})$ is an inverse sequence.
Clearly (see (d) and (e) of Lemma~\ref{L3.1}), $\lim\bold
P_x'=\pi^{-1}(x)$.  Now $\bold P_x'$ is an inverse sequence of
contractible polyhedra.  Hence (II) is true.

To get at (III), first observe that by (f) of Lemma~\ref{L3.1}, the fiber
$(\pi\vert_{A_{k}})^{-1}(x)$ is the limit of the inverse sequence
$(S_{x,i},g_i^{i+1})$.
On the other hand, for each $i\in\N$, $B_{x,i}\subset
P_{x,i}\subset B_{x,i}^{\#}$, $g_i^{i+1}(P_{i+1}^{(l_k)})\subset
P_i^{(l_k)}$, and $g_i^{i+1}(B_{x,i+1}^{\#})\subset B_{x,i}$.  So
one deduces that
$$g_i^{i+1}(P_{x,i+1}^{(l_k)})\subset
g_i^{i+1}(B_{x,i+1}^{\#})\cap g_i^{i+1}(P_{i+1}^{(l_k)})\subset
B_{x,i}\cap P_i^{(l_k)}\subset P_{x,i}\cap P_i^{(l_k)}=P_{x,i}^{(l_k)}.$$
Thus $\bold P_x^{'(l_k)}:=(P_{x,i}^{(l_k)},g_i^{i+1})$ is an inverse
sequence of compact polyhedra.  Since $S_{x,i}\subset
P_{x,i}^{(l_k)}$ and
$$g_i^{i+1}(P_{x,i+1}^{(l_k)})\subset B_{x,i}\cap P_i^{(l_k)}=S_{x,i},$$
it is clear that $\lim\bold P_x^{'(l_k)}$ is the same as the limit of
the inverse sequence $(S_{x,i},g_i^{i+1})$, i.e., that
$$(\pi|_{A_k})^{-1}(x)=\lim \ (S_{x,i},g_i^{i+1})=\lim \ (P_{x,i}^{(l_k)},g_i^{i+1}).$$

We shall show that for each $i\in\N$, if
$0\leq r\leq l_k-1$ and $h:S^r\ra P_{x,i}^{(l_k)}$ is a map, then $h$
is homotopic to a constant map.
Since $\dim S^r=r<l_k$, $h$ is homotopic in $P_{x,i}^{(l_k)}$ to a map
that carries $S^r$ into $P_{x,i}^{(l_k-1)}$ (see remark about stability theory).  But $P_{x,i}$ is
contractible, so the inclusion $P_{x,i}^{(l_k-1)}\hookrightarrow
P_{x,i}^{(l_k)}$ is null-homotopic.  This shows that $h:S^r\ra
P_{x,i}^{(l_k)}$ is null-homotopic.  So all fibers of $\pi\vert_{A_k}$
are $\UV^{l_k-1}$.

\vspace{2mm}

To prove (IV), we need to show that any fiber of
$\pi\vert_{A_k}$ is $\Z/p$ -acyclic, i.e., for infinitely many
indexes $i$, the map $g_i^{i+1}\vert_{P_{x,i+1}^{(l_k)}}:
P_{x,i+1}^{(l_k)} \ra P_{x,i}^{(l_k)}$ induces the zero-homomorphism
of cohomology groups $\mathrm{H}^m (P_{x,i}^{(l_k)};\Z/p) \ra
\mathrm{H}^m (P_{x,i+1}^{(l_k)};\Z/p)$, for all $m \in \N$ (we need
not worry about $m=0$ because the $P_{x,i}^{(l_k)}$'s are
$(l_k-1)$-connected, so their reduced zero-cohomology groups are
trivial). We will be focusing on those indexes $i$ for which $g_i^{i+1}\vert_{P_{x,i+1}} \simeq \omega \circ \bar{f_i}\vert_{P_{x,i+1}}$, as mentioned in the conditions of Corollary~\ref{C2.2}.

It is, in fact, enough to show that the map $g_i^{i+1}\vert_{P_{x,i+1}^{(l_k)}}:
P_{x,i+1}^{(l_k)} \ra P_{x,i}^{(l_k)}$ induces the zero-homomorphism of homology groups with $\Z/p$-coefficients. Here is why this is true. Notice that each of
 $P_{x,i+1}$ and $P_{x,i}$ is a closed vertex star (in the coarser triangulation), subdivided so that each original simplex of the vertex star is collapsible as a simplicial complex. So Lemma~\ref{H_n(EW)} (for $n=l_k$) is true for both $|M|=P_{x,i+1}$ and $|M|=P_{x,i}$. Therefore
 property (I) of  Lemma~\ref{H_n(EW)} is true for both $P_{x,i+1}^{(l_k)}$ and $P_{x,i}^{(l_k)}$, and both are $(l_k-1)$-connected. Therefore by the Universal Coefficients Theorem for homology and cohomology we have
$$H_{m}(P_{x,i+1}^{(l_k)};\Z/p)\cong H_{m}(P_{x,i+1}^{(l_k)})\otimes \Z/p,\ \  \forall m\geq 1,
\text{ and}$$
$$ H^{m}(P_{x,i+1}^{(l_k)};\Z/p)\cong \Hom (H_{m}(P_{x,i+1}^{(l_k)}), \Z/p),\ \ \forall m\geq 1,$$ 
and for $P_{x,i}^{(l_k)}$ analogously, and these expressions are non-zero only for $m=l_k$.
So if $g_i^{i+1}\vert_{P_{x,i+1}^{(l_k)}}:
P_{x,i+1}^{(l_k)} \ra P_{x,i}^{(l_k)}$ induces the zero-homomorphism  $H_{l_k}(g_i^{i+1};\Z/p): H_{l_k}(P_{x,i+1}^{(l_k)};\Z/p) \ra H_{l_k}(P_{x,i}^{(l_k)};\Z/p)$, then for any $\vp \in \Hom (H_{l_k}(P_{x,i}^{(l_k)}),\Z/p)$, we have $\vp\circ H_{l_k}(g_i^{i+1})=0 \in \Hom (H_{l_k}(P_{x,i+1}^{(l_k)}),\Z/p)$, that is, the induced homomorphism 
$ H^{l_k}(g_i^{i+1};\Z/p): H^{l_k}(P_{x,i}^{(l_k)};\Z/p) \ra H^{l_k}(P_{x,i+1}^{(l_k)};\Z/p)$ is the zero-homomorphism.

\vspace{2mm}

So let us show that $H_{l_k}(g_i^{i+1};\Z/p): H_{l_k}(P_{x,i+1}^{(l_k)};\Z/p) \ra H_{l_k}(P_{x,i}^{(l_k)};\Z/p)$ is the zero-ho\-mo\-mor\-phism.
Before proceeding, note that by Remark~\ref{R-EW},
given an EW-resolution $\omega :\EW(P_i,\Z/p,l_k) \ra P_i$, we know that
$\omega^{-1}(P_{x,i})=\EW(P_{x,i},\Z/p,l_k)$, so \\
$\omega\vert_{\omega^{-1}(P_{x,i})}:\EW(P_{x,i},\Z/p,l_k)\ra P_{x,i}$ is also an EW-resolution.

\vspace{2mm}

Let $\sa$ be any $(l_k+1)$-simplex of $P_{x,i+1}$, and let $g_\sa$ be a generator of $H_{l_k}(\partial\sa)$. Let $\mu: \partial \sa \hookrightarrow
P_{x,i+1}^{(l_k)}$, $\nu : P_{x,i+1}^{(l_k)} \hookrightarrow
P_{x,i+1}$, and $\la : P_{x,i}^{(l_k)} \hookrightarrow
\EW(P_{x,i},\Z/p,l_k)$ be the inclusions. Notice that $\omega
\circ \bar{f_i} (P_{x,i+1}) \subset P_{x,i}$ implies that
$\bar{f_i} (P_{x,i+1}) \subset \EW(P_{x,i},\Z/p,l_k)$, and
since $\bar{f_i}$ is a cellular map, we also have $\bar{f_i}
(P_{x,i+1}^{(l_k)}) \subset \EW(P_{x,i},\Z/p,l_k)^{(l_k)} =
P_{x,i}^{(l_k)}$.
{\small
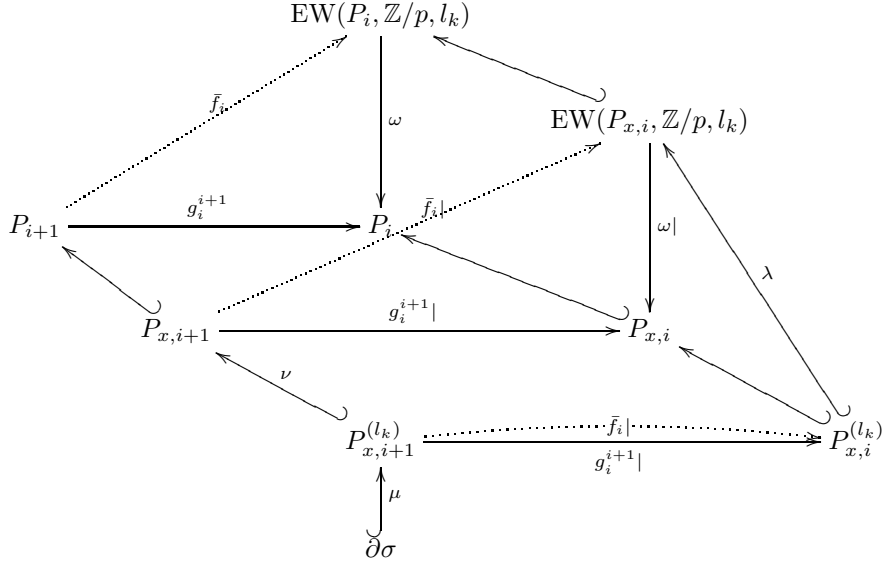
\begin{displaymath}
\xymatrix{
&& \EW(P_i,\Z/p,l_k) \ar[dd]^{\omega} &&\\
&&& \EW(P_{x,i},\Z/p,l_k) \ar[dd]^{\omega\vert} \ar@{_{(}->}[ul]  &\\
P_{i+1}\ar@{..>}[rruu]^{\qquad \qquad\bar f_i} \ar[rr]^{g_i^{i+1}}&& P_i &&\\
&P_{x,i+1} \ar@{_{(}->}[ul] \ar[rr]^{g_i^{i+1}\vert} \ar@{..>}[rruu]^{\qquad \qquad\bar f_i\vert} && P_{x,i}\ar@{_{(}->}[ul] &\\
&&P_{x,i+1}^{(l_k)} \ar@{_{(}->}[ul]_{\nu} \ar[rr]_{g_i^{i+1}\vert} \ar@/^/@{..>}[rr]|{\bar f_i\vert} && P_{x,i}^{(l_k)} \ar@{_{(}->}[ul] \ar@{_{(}->}[uuul]_{\la} \\
&& \partial \sa \ar@{^{(}->}[u]_{\mu}&&
}
\end{displaymath}
}
Since  Lemma~\ref{algebra} is true for $|M|=P_{x,i}$ and $n=l_k$, we have 
$\bar{f_i}\vert_{P_{x,i+1}}\circ \nu \circ \mu =
\la \circ \bar{f_i}\vert_{P_{x,i+1}^{(l_k)}}\circ \mu$ is null-homotopic, and 
$$H_{l_k}(\bar{f_i}\vert_{P_{x,i+1}^{(l_k)}}\circ
\mu)(g_\sa)=\sum_{e=1}^r \e_eh_e \ \in \ H_{l_k}(P_{x,i}^{(l_k)}),
\text{ where }
\e_e \equiv 0 \ (\mod \ p).$$

By Lemma~\ref{H_n(EW)} applied to $P_{x,i+1}$ with $n=l_k$, we can select $\sa_1, \ldots , \sa_s$ to be some $(l_k+1)$-simplexes of $P_{x,i+1}$ so that the images $g_1, \ldots , g_s$ of the generators of $H_{l_k}(\partial\sa_1), \ldots , H_{l_k}(\partial\sa_s)$ induced by the inclusions $\partial\sa_j\hookrightarrow P_{x,i+1}^{(l_k)}$ form a basis for $H_{l_k}( P_{x,i+1}^{(l_k)})$.

Then for any $g\in H_{l_k}(P_{x,i+1}^{(l_k)})$,
$$H_{l_k}(\bar{f_i}\vert_{P_{x,i+1}^{(l_k)}})( g)=H_{l_k}(\bar{f_i}\vert_{P_{x,i+1}^{(l_k)}}) (\sum_{j=1}^s m_j
g_{j})= \sum_{j=1}^s m_j(\sum_{e=1}^r \e_{j,e}h_e),$$
where $m_j\in\Z$, and $\e_{j,e} \equiv 0 \ (\mod \ p)$. 

\vspace{2mm}

Finally, since we know that $g_i^{i+1}\vert_{P_{x,i+1}} \simeq \omega \circ \bar{f_i}\vert_{P_{x,i+1}}$ and
$\omega \vert_{P_{x,i}^{(l_k)}} = id$, we have that  $g_i^{i+1}\vert_{P_{x,i+1}^{(l_k)}} \simeq \bar{f_i}\vert_{P_{x,i+1}^{(l_k)}} $. Therefore $H_{l_k}(g_i^{i+1}\vert_{P_{x,i+1}^{(l_k)}} )=H_{l_k}(\bar f_i\vert_{P_{x,i+1}^{(l_k)}} )$, so the last equation implies that $H_{l_k}(g_i^{i+1}\vert_{P_{x,i+1}^{(l_k)}};\Z/p )$ is the zero-homomorphism. $\square$

\section{Proof of Theorem~\ref{T1.1}}

\noindent \textsl{Proof of Theorem~\ref{T1.1}}: \noindent Choose
a function $\nu:\N\ra\N$ such that for each $i\in\N$,
\begin{enumerate}
\item[(i)] $\nu(i)\leq i$, and
\item[(ii)]  $\nu^{-1}(i)$ is infinite.
\end{enumerate}
One may assume that $X\subset Q=$ Hilbert cube. We are
going to prove the existence for each $k\in\Bbb N\cup\{\infty\}$
of a certain sequence $\mathcal{S}_j=( n_j,\ (P^k_j),\
\varepsilon_j, \ \delta_j,\ (\widetilde\tau^k_j)), \ (\tau^k_j)),\ j\in\N$, of entities,
and a sequence of maps $(g_j^{j+1}),\ j\in\N$, such that:
\begin{itemize}
\item $n_j\in\N$;
\item $P_j^1\subset P_j^2\subset \dots
 \subset P_j^\infty$ are compact subpolyhedra of $I^{n_j}$;
\item  $\varepsilon_j,\ \delta_j>0  $;
\item $\widetilde\tau^\infty_j$ is a triangulation of
$ P^\infty_j$ , $\tau^\infty_j$  is a subdivision of $\widetilde\tau^\infty_j$,\\
$ \widetilde\tau^k_j=\widetilde\tau^\infty_j|_{P^k_j}$  is a triangulation of
 $ P^k_j$ , $\tau^k_j=\tau^\infty_j|_{P^k_j}$ is a subdivision of $\widetilde\tau^k_j$,\\
(we will consider $P^\infty_j=(P^\infty_j, \tau^\infty_j)$  and 
$P^k_j=(P^\infty_j, \tau^k_j)$);
\item $\ g^{j+1}_j\colon P^\infty_{j+1}\rightarrow P^\infty_j$
  is a simplicial
map relative to $ \tau_{j+1}^\infty$ and $\tau_{j}^\infty$.
\end{itemize}

A diagram that might help:
{\footnotesize
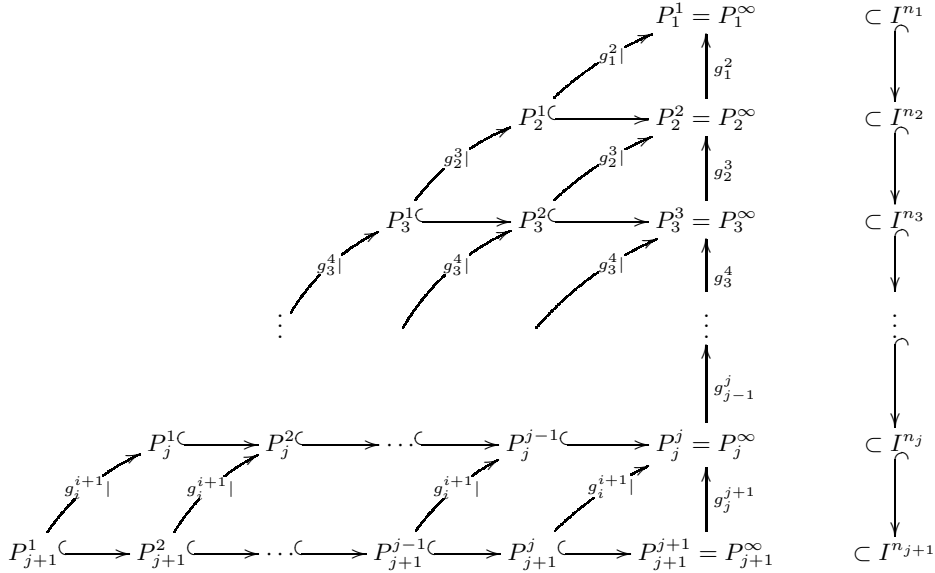
\begin{displaymath}
\xymatrix{
&&&&&                                                                                            P_1^1= P_1^\infty  &\subset I^{n_1}\ar@{^{(}->}[d]\\
&&&&                           P_2^1 \ar@/^/[ur]|{g^2_1\vert} \ar@{^{(}->}[r] & P_2^2=P_2^\infty \ar[u]_{g^2_1}  &\subset I^{n_2} \ar@{^{(}->}[d]\\
&&&    P_3^1 \ar@/^/[ur]|{g^3_2\vert}  \ar@{^{(}->}[r] & P_3^2 \ar@/^/[ur]|{g^3_2\vert}  \ar@{^{(}->}[r] & P_3^3=P_3^\infty \ar[u]_{g^3_2} &\subset {I^{n_3}}\ar@{^{(}->}[d]\\
&&\vdots\ar@/^/[ur]|{g^4_3\vert}&\ar@/^/[ur]|{g^4_3\vert}&\ar@/^/[ur]|{g^4_3\vert}&\vdots \ar[u]_{g^4_{3}}&\vdots\ar@{^{(}->}[d]\\
&P_j^1  \ar@{^{(}->}[r] & P_j^2  \ar@{^{(}->}[r] & \cdots \ar@{^{(}->}[r]  &P_j^{j-1} \ar@{^{(}->}[r] & P_j^j=P_j^\infty \ar[u]_{g^j_{j-1}}  &\subset {I^{n_j}}\ar@{^{(}->}[d]\\
 P_{j+1}^1 \ar@{^{(}->}[r] \ar@/^/[ur]|{g^{i+1}_i\vert} & P_{j+1}^2 \ar@{^{(}->}[r] \ar@/^/[ur]|{g^{i+1}_i\vert} &\cdots \ar@{^{(}->}[r] & P_{j+1}^{j-1} \ar@{^{(}->}[r]\ar@/^/[ur]|{g^{i+1}_i\vert}& P_{j+1}^{j} \ar@{^{(}->}[r] \ar@/^/[ur]|{g^{i+1}_i\vert}&  P_{j+1}^{j+1}=P_{j+1}^\infty
  \ar[u]_{g^{j+1}_{j}}& \subset {I^{n_{j+1}}}  \\
}
\end{displaymath}
}

We shall require that for each $j\in \N$ and $k\in\N$:
\begin{enumerate}
\item[$(1)_{j>1}$] $n_{j-1}<n_j$;
\item[$(2)_{j\ge 1}$] if $j\leq k<\infty$, then $P_j^k=P_j^\infty$
 and $P_j^r\subset\int_{I^{n_j}}P_j^{r+1}$ whenever $r<j$;
\item[$(3)_{j\ge 1}$] $X\subset\int_Q(P_j^\infty\times Q_{n_j})
\subset N(X,\frac{2}{j})$, and,\\ whenever $k<j$,
 $X_k\subset\int_Q(P_j^k\times Q_{n_j})\subset
N(X_k,\frac{2}{j})$;
\item[$(4)_{j>1}$]
$p_{n_{j-1}}(P_{j}^k)\subset\int_{I^{n_{j-1}}}P_{j-1}^k$;
\item[$(5)_{j> 1}$] if  $u$, $v\in Q$ and
$\rho(u,v)\leq\e_{j}$, then
$\rho(p_{n_{j-1}}(u),p_{n_{j-1}}(v))<\da_{j-1}$;
\item[$(6)_{j\ge 1}$]  $\frac{9}{2^{n_j}}<\e_j$;
\item[ $(7)_{j\ge 1}$] $\da_j<\frac{1}{2^{n_j-1}}$;
\item[ $(8)_{j\ge 1}$] $\tau^\infty_j|_{|\ga|}$ is collapsible $\forall\ga\in \widetilde\tau^\infty_j$ and 
$\mesh\tau_j^\infty<\frac{\da_j}{2}$;
\item[$(9)_{j\ge 1}$]
if $x\in X$, then there exists a contractible subpolyhedron $P_{x,j}^\infty$ of
$P_j^\infty$, which is the closed star of a vertex in the triangulation $\widetilde\tau^\infty_j$, i.e., $P_{x,j}^\infty =\overline\St(v, \widetilde\tau^\infty_j)$ for some $v\in (\widetilde\tau^\infty_j)^{(0)}$, and such that
$\overline N(p_{n_j}(x),2\da_j)\cap P_j^\infty \subset P_{x,j}^\infty\subset\overline
N(p_{n_j}(x),\e_j)\cap P_j^\infty$;\\
($P^\infty_{x,j}$ is considered with the triangulation $\tau^\infty_j$, so it is a subdivided vertex star);\\
if $k< j$, and
$x\in X_k$, then there exists a contractible subpolyhedron $P_{x,j}^k$ of
$P_j^k$, which is the closed star of a vertex in the triangulation $\widetilde\tau^k_j$, i.e., $P_{x,j}^k =\overline\St(v, \widetilde\tau^k_j)$ for some $v\in (\widetilde\tau^k_j)^{(0)}$, and such that
$\overline N(p_{n_j}(x),2\da_j)\cap P_j^k\subset P_{x,j}^k\subset\overline
N(p_{n_j}(x),\e_j)\cap P_j^k$;\\
($P^k_{x,j}$ is considered with the triangulation $\tau^k_j$).
This statement is also true when $k\geq j$, because then $P_{x,j}^k=P_{x,j}^\infty$, $P_j^k=P_j^\infty$ and $X_k\subset X$;
\item[$(10)_{j>1}$] whenever $x\in P_{j}^\infty$ there exists a
simplex $\sa$ of $\tau_{j-1}^\infty$ such that
$g^j_{j-1}(x)\in\sa$, and $p_{n_{j-1}}(x)$ lies in
$N(\sa,\frac{\da_{j-1}}{2})$ (and therefore, it follows from here
and $(8)_{j-1}$ that, $\rho(g_{j-1}^j(x),p_{n_{j-1}}(x))<
\da_{j-1}/2+\da_{j-1}/2=\da_{j-1}$ for all $x\in P_j^\infty$);
\item[$(11)_{j>1}$] $g^{j}_{j-1}(P_{j}^k)\subset P_{j-1}^k$; and
\item[$(12)_{j>1}$] $g_{j-1}^{j}\vert_{P_{j}^{\nu(j-1)}} \simeq
\omega \circ \bar{f}_{j-1}$, where $\omega :
\EW(P_{j-1}^{\nu(j-1)},\Z/p,l_{\nu(j-1)}) \ra
P_{j-1}^{\nu(j-1)}$ is an Edwards--Walsh projection, and $
\bar{f}_{j-1} : P_{j}^{\nu(j-1)} \ra
\EW(P_{j-1}^{\nu(j-1)},\Z/p,l_{\nu(j-1)})$ is a cellular
map. Moreover, for all $x \in X_{\nu(j-1)}$, we have that \
$\omega \circ \bar{f}_{j-1}(P_{x,j}^{\nu(j-1)} ) \subset
P_{x,j-1}^{\nu(j-1)} $, and $g_{j-1}^{j}\vert_{P_{x,j}^{\nu(j-1)}} \simeq \omega \circ
\bar{f}_{j-1}\vert_{P_{x,j}^{\nu(j-1)}}$.
\end{enumerate}

{\small
\begin{displaymath}
\xymatrix{
&&  \EW(P_{j-1}^{\nu(j-1)},\Z/p,l_{\nu(j-1)}) \ar[dd]^{\omega} &&\\
&&&  \EW(P_{x,j-1}^{\nu(j-1)},\Z/p,l_{\nu(j-1)}) \ar[dd]^{\omega\vert} \ar@{_{(}->}[ul]  &\\
P_{j}^{\nu(j-1)}\ar@{..>}[rruu]^{\qquad \qquad\bar f_{j-1}} \ar[rr]^{g_{j-1}^{j}|}&& P_{j-1}^{\nu(j-1)} &&\\
& P_{x,j}^{\nu(j-1)}\ar@{_{(}->}[ul] \ar[rr]^{g_{j-1}^{j}\vert} \ar@{..>}[rruu]^{\qquad \qquad \bar{f}_{j-1}\vert} &&  P_{x,j-1}^{\nu(j-1)}\ar@{_{(}->}[ul] &
}
\end{displaymath}
}

Before proving the existence of such data, let us see
why they would imply the conclusion of Theorem~\ref{T1.1}.  For each
$i\in\N$, let $P_i^\infty$ correspond to $P_i$ from the statement of Lemma~\ref{L3.1}. Applying (5), (1), (6), (10) and (7), one sees that the conditions (i)--(v) of Lemma~\ref{L3.1} are clearly true. Condition $(4)_{i+1}$ implies (vi) and one
may use (3) to see that
 $$X=\bigcap_{i=1}^{\infty}
P_i^\infty \times Q_{n_i}.$$
Let $$Z:=\lim(P_i^\infty,g^{i+1}_{i}).$$
Surely $Z$ is
a metrizable compactum, and we get the map $\pi:Z\ra X$ defined by
the formula given in Lemma~\ref{L3.1} (b).

To see that $\pi$ is surjective, for each $i\in\N$ let
$T_i=P_i=P_i^\infty$ (in Lemma~\ref{L3.1}). According to the notation of the last part
of Lemma~\ref{L3.1}, one sees that for $x\in X$,
$S_{x,i}=B_{x,i}=\overline N(p_{n_i}(x),2\da_i)\cap P_i^\infty$ (while $B_{x,i}^\#=\overline N(p_{n_i}(x),\e_i)\cap P_i^\infty$).  Notice that the first part of $(2)_i$ together with
$(3)_i$  implies

\begin{enumerate}
\item[(13)] \ $p_{n_i}(X)\subset \int_{I^{n_i}}P_i^\infty$, and  $\forall k\in \N$,
              $p_{n_i}(X_k)\subset \int_{I^{n_i}}P_i^k$.
\end{enumerate}

So $p_{n_i}(x)\in P_i^\infty$ and therefore
$p_{n_i}(x)\in B_{x,i}$, showing that the latter is not empty. The
map $\tilde\pi$ is the same as $\pi$ in this setting, so (g) of
Lemma~\ref{L3.1} shows that $\pi$ is surjective.

One then checks that all the hypotheses of Corollary~\ref{C2.2} except
for the very last one (which we do not need yet) are also
satisfied. Thus (I)--(III) of Corollary~\ref{C2.2} hold true, so $\pi$ is a cell-like map,
and we are assured of the existence of the closed subspaces $A_k$,
$k\geq 1$, where
$$A_k:=\lim((P_i^\infty)^{(l_k)},g^{i+1}_{i}),$$ as required by
Theorem~\ref{T1.1} so that $\dim A_k\leq l_k$, and when $k\in \N$, $\pi$
carries $A_k$ in a $\UV^{l_k-1}$ manner onto $X$.

\vspace{2mm}

We must identify the closed subspaces $Z_1\subset Z_2\subset
\dots$ of $Z$, prove they satisfy (a)--(c) of Theorem~\ref{T1.1}, and
show that $Z_k\subset A_k$ when $k\in\N$.  Fix $k\in\N$. In the
last part of Lemma~\ref{L3.1}, instead of putting $T_i=P_i^\infty$, as we
just did to obtain $Z$, $\pi$, and the sets $A_k$, this time put
$T_i=(P_i^k)^{(l_k)}$.  Using (11), the fact that
$\tau_i^k=\tau_i^\infty|_{P_i^k}$, and that $g_i^{i+1}$ is simplicial
from $\tau_{i+1}^\infty$ to $\tau_i^\infty$, one sees that
\begin{enumerate}
\item[(14)] $g_i^{i+1}((P_{i+1}^k)^{(l_k)})\subset(P_i^k)^{(l_k)}$.
\end{enumerate}
Now let
$$Z_k:=\lim((P_i^k)^{(l_k)},g_i^{i+1}),$$
i.e., $\bold T_k=((P_i^k)^{(l_k)},g_i^{i+1})$, and $Z_k=\lim \bold T_k$.
Using (2) we see that $P_i^{k}\subset
P_i^\infty$ for all $i\in\N$. Of course,
$(P_i^k)^{(l_k)}\subset(P_i^\infty)^{(l_k)}$, and we deduce that
$Z_k\subset A_k$ as requested in Theorem~\ref{T1.1}.  Moreover, $\dim
Z_k\leq\dim A_k\leq l_k$, so (a) of Theorem~\ref{T1.1} has been resolved.
It is also clear that $Z_1\subset Z_2\subset\dots$ as required by Theorem~\ref{T1.1}.

Next put $\tilde\pi_k=\pi|Z_k:Z_k\ra X$.  If $(a_1,a_2,\dots)$ is
a thread of $Z_k$, then $a_i\in P_i^k$ for each $i\in\N$.  Taking
into account (b) of Lemma~\ref{L3.1}, as well as (3)$_i$ which implies that
\begin{enumerate}
\item[(15)] $\displaystyle X_k=\bigcap_{i=1}^{\infty} P_i^k \times Q_{n_i},$
\end{enumerate}
one sees that
$\tilde\pi_k(Z_k)\subset X_k$.

Suppose now that $x\in X_k$.  With the choice of
$T_i=(P_i^k)^{(l_k)}$, the sets $S_{x,i}$ in the last part of Lemma~\ref{L3.1} become
$S_{x,i}=B_{x,i}\cap(P_i^k)^{(l_k)}.$

If we can show that for each $i\in\N$, $S_{x,i}\neq\emptyset$,
then (g) of Lemma~\ref{L3.1} would yield  $\tilde\pi_k(Z_k)\supset X_k$. Indeed, it
is sufficient to show that
$B_{x,i}^{\#}\cap(P_{i}^k)^{(l_k)}\neq\emptyset$, since $B_{x,i+1}^\#\cap(P_{i+1}^k)^{(l_k)}$
maps into $S_{x,i}$ under $g_i^{i+1}$(see (c) of Lemma~\ref{L3.1} and (14)).
Because of (15), $x\in P_i^k\times Q_{n_i}$, so $p_{n_i}(x)\in
P_i^k$. Applying (6)$_i$--(8)$_i$, we find a vertex
$v\in(P_i^k)^{(0)}\subset (P_i^k)^{(l_k)}$ such that
$\rho(p_{n_i}(x),v)<\frac{\da_i}{2}<\frac{1}{2^{n_i}}<\e_i.$
This means $v\in
B_{x,i}^{\#}\cap(P_i^k)^{(l_k)}$, i.e.,
$B_{x,i}^{\#}\cap(P_i^k)^{(l_k)}\neq\emptyset$. Therefore (b) of Theorem~\ref{T1.1} is true.

\vspace{2mm}

Finally, after replacing $A_k$ from the statement of Corollary~\ref{C2.2} with $Z_k$, the ultimate condition of Corollary~\ref{C2.2}, involving infinitely many indexes, is now operative because
of (i) and (ii) of this section, and $(12)$ for $\nu(i-1)=k$.  If we
apply (IV) of Corollary~\ref{C2.2}, then we find that
$\tilde\pi_k=\pi\vert_{Z_k}:Z_k\ra X_k$ is a $\Z/p$-acyclic map.
Thus, our proof of Theorem~\ref{T1.1} will be complete once we have obtained
the information in statements (1)--(12).

\vspace{4mm}

\noindent \textbf{Inductive construction begins}: For the basis of the induction ($j=1$), we choose $n_1=l_1$ and $P_1^k=I^{n_1}=I^{l_1}$ for all $k\in\N\cup\{\infty\}$. Thus $(2)_1$ and $(3)_1$ are satisfied.
Next choose any $\e_1 > \frac{9}{2^{l_1}}$, so $(6)_1$ is satisfied. It remains to produce $\da_1 >0$ and triangulations $\widetilde\tau_1^\infty$ and $\tau_1^\infty$  of $P_1^\infty=I^{l_1}$ so that $(7)_1$--$(9)_1$ are satisfied.

Begin by taking a triangulation $\widetilde\tau_1^\infty$ of
$P_1^\infty$ such that $\mesh \widetilde\tau_1^\infty< \frac{\e_{1}}{2}$. The open stars of the vertices in $\widetilde\tau_1^\infty$ form a cover for $P_{1}^\infty=I^{l_1}$. Note that these open stars are truly open sets in $I^{l_1}$.
For any $x \in X$,
there exists a vertex $v$ of $\widetilde\tau_1^\infty$  such that
$p_{l_{1}}(x) \in \St(v, \widetilde\tau_1^\infty)$.
Note
that for any $y \in \St(v, \widetilde\tau_1^\infty)$, $\rho (y,p_{l_{1}}(x))\leq
2\mesh \widetilde\tau_1^\infty<\e_{1}$, so $\overline{\St}(v, \widetilde\tau_1^\infty) \subset
\overline{N}(p_{l_{1}}(x), \e_{1})= \overline{N}(p_{l_{1}}(x), \e_{1})\cap P_1^\infty$.

Since  $\scu := \{\St (v, \widetilde\tau_1^\infty) \vert \ v \in (\widetilde\tau_1^\infty)^{(0)} \}$ is a cover for $P_1^\infty$ which is compact, let $\la$ be a Lebesgue number of $\scu$.
Pick a $\da_{1}>0$ such that
$4\da_{1}\ < \min \left\{ \la, \frac{4}{2^{l_{1}-1}}\right\}.$
Now $(7)_{1}$ is also satisfied. Then for
any $x \in X$, the closed ball
$\overline{N}(p_{l_1}(x),2\da_{1})$ is contained in some $\St(v, \widetilde\tau_1^\infty)$, for a vertex $v\in (\widetilde\tau_1^\infty)^{(0)}$. Pick one such star, and call its
closure $P_{x,1}^\infty$. Notice that $P_{x,1}^\infty$ is contractible.
Thus we get $(9)_{1}$ for $x\in X$:
$\overline N(p_{l_1}(x),2\da_1)= \overline N(p_{l_1}(x),2\da_1) \cap P_1^\infty \subset P_{x,1}^\infty\subset\overline
N(p_{l_1}(x),\e_1)\cap P_1^\infty$.
Finally,
choose a triangulation $\tau_{1}^\infty$ so that it refines
$\widetilde\tau_1^\infty$, and so that $(8)_{1}$ is satisfied.

\vspace{2mm}

 Assume that we have completed the construction of
$\mathcal{S}_j$ for $1\leq j \leq i$, and $g_j^{j+1}$ for $1\leq j
\leq i-1$. Choose an open cover $\mathcal{V}$ of $P_i^\infty$
having the property that $\mesh\mathcal{V}<\frac{\da_i}{2}$.  Then
select a finer open cover $\mathcal{W}$ such that any two
$\mathcal{W}$-near maps of any space into $P_i^\infty$ are
$\mathcal{V}$-homotopic. Let $\tau$ be a subdivision of
$\tau_i^\infty$ such that $N(\overline{\St}(v, \tau), \tilde \e)$
lies in an element of $\mathcal{W}$, for every vertex $v \in
\tau^{(0)}$, where $\tilde \e >0 $ is chosen so that: for any
principal simplex $\sa$ of the triangulation $\tau$, all of the points
of the open neighborhood $N(\sa, \tilde \e)$  are at most one
(principal) simplex away from $\sa$ (i.e., if $u \in N(\sa,\tilde
\e) \setminus \sa$, then $u \in \gamma =$ a neighboring principal
simplex of $\sa$). (Surely this $\tilde \e$ exists because
$P_i^{\infty}$ is compact. Also, it is clear that $\tilde \e \leq
\mesh \tau$, and that it would be enough to choose $\tau$ so that
$2(\mesh \tau + \tilde \e) <$ some fixed Lebesgue number of
$\mathcal{W}$. Also note that $\tau$ can be chosen so that $\tau|_{|\ga|}$ is still  collapsible, $\forall \ga \in \widetilde\tau_i^\infty$.)

If $i=1$, replace $\tau_1^\infty$ by $\tau$, but continue to use the notation $\tau_1^\infty$ for it.
Note that properties $(8)_1$ and $(9)_1$, which are the only ones affected by this change, are still true.

 If $i>1$, choose a map $\mu:P_i^\infty\ra P_i^\infty$ which is
simplicial from $\tau$ to $\tau_i^\infty$ and which is a
simplicial approximation to the identity on $P_i^\infty$.  Then
the map $g_{i-1}^i \circ \mu$ is simplicial from $\tau$ to
$\tau_{i-1}^\infty$, and $\bar{f}_{i-1} \circ \mu\vert_{P_i^{\nu(i-1)}}$ is cellular
with respect to the triangulation on $P_i^{\nu(i-1)}$ induced by
$\tau$ for $\bar f_{i-1}$. If
we replace $g_{i-1}^i$ by $g_{i-1}^i \circ \mu$,  $\bar{f}_{i-1}$
by $\bar{f}_{i-1} \circ \mu\vert_{P_i^{\nu(i-1)}}$, and
$\tau_i^\infty$ by $\tau$, then all the conditions (1)--(12) for
index $i$ still prevail (the only ones affected being
(8)$_i$--(12)$_i$). So we assume that these replacements have been
made, but continue to use $g_{i-1}^i$, $\bar{f}_{i-1}$ and
$\tau_i^\infty$ to denote the respective bonding map, cellular map
in (12)$_i$ and triangulation.

\vspace{4mm}

\noindent \textbf{Construction of the polyhedra $P_{i+1}^k$
and the bonding map $g_i^{i+1}$ begins}. Apply Lemma~\ref{L8.2}
to $X_{\nu(i)}$, which has $\dim _{\Z/p} X_{\nu(i)} \leq
l_{\nu(i)}$, where $\nu(i)\leq i$, and (using (13)) the map $p_{n_i}|_{X_{\nu(i)}}\ : X_{\nu(i)} \to
P_i^{\nu(i)}$, to produce a map $f':X_{\nu(i)} \to
\EW(P_i^{\nu(i)},\Z/p ,l_{\nu(i)})$ such that for any $x \in
X_{\nu(i)}$, when $p_{n_i}(x)$ lies in a particular simplex of $P_i^{\nu(i)}$, then so does $\omega \circ f' (x)$. There is a principal simplex $\sa_x$ of $P_i^{\nu(i)}$ that contains both
$\omega\circ f'(x)$ and $p_{n_i}(x)$. We can extend $f'$ over an open
neighborhood $\widetilde U$ of $X_{\nu(i)}$ in the Hilbert cube
$Q$, to get a map $f'' : \widetilde U \to
\EW(P_i^{\nu(i)},\Z/p ,l_{\nu(i)})$.
\begin{displaymath}
\xymatrix{
&& \EW (P_i^{\nu(i)},\Z/p ,l_{\nu(i)}) \ar[dd]_{\omega}\\
&&\\
\widetilde U \ar@{..>}[rruu]|{f''} & X_{\nu(i)} \ar@{_{(}->}[l] \ar[r]_{p_{n_{i}\vert}} \ar@{-->}[ruu]_{\!f'}& P_i^{\nu(i)}\\
}
\end{displaymath}
Now we can find a
neighborhood $U$ of $X_{\nu(i)}$ in $\widetilde U$ such that:
\begin{enumerate}
\item[(16)] for any $u \in U$, $\omega\circ f'' (u)$ and $p_{n_i}(u)$
belong to the open $\tilde \e$--neighborhood of some principal
simplex $\sa_x$ of $P_i^{\nu(i)}$.
\end{enumerate}
Here is how we find $U$: since $p_{n_i}$ is continuous (on $Q
\supset \widetilde U$), for any $x\in X_{\nu(i)}$, and for the
above $\tilde \e$, there exists an open neighborhood $\widetilde
Q_x$ of $x$ in $\widetilde U$ such that $p_{n_i}(\widetilde Q_x)
\subset N(\sa_x, \tilde \e)$.
Since $\omega\circ f'(x)\in\sa_x$, then $f'(x)\in\omega^{-1}(\sa_x)\subset\omega^{-1}(N(\sa_x,\tilde\e))$.
Now $f''(x)=f'(x)$, so the continuity of $f''$ guarantees an open neighborhood
$\bar Q_x$ of $x$ with
$f''(\bar Q_x)\subset\omega^{-1}(N(\sa_x,\tilde\e))$.  Of
course, $\omega\circ f''(\bar Q_x)\subset N(\sa_x,\tilde\e)$.

Now let
$Q_x:=\widetilde Q_x \cap \bar Q_x$ and  define $U:= \bigcup_{x
\in X} Q_x$. Clearly this $U$ has the needed property.

\vspace{4mm}

Using the uniform continuity of $p_{n_i}$ on $Q$, choose $\e_{i+1}$ so
that (5$)_{i+1}$ holds: if $u,v \in Q$ are such that $\rho (u,v) <
\e_{i+1}$, then $\rho (p_{n_i}(u),p_{n_i}(v))\ < \da_i$.

\vspace{4mm}

In order to choose $n_{i+1}$: notice that one may find $m_0\in\N$
such that if $m\geq m_0$, then $X\subset p_m(X)\times Q_m\subset
N(X,\frac{2}{i+1})$, and for all $k\leq i$, $X_k\subset
p_m(X_k)\times Q_m \subset N(X_k,\frac{2}{i+1})$. Define $n_{i+1}
> \max \ \{l_{i+1}-1,n_i, m_0, \log_2 (\frac{9}{\e_{i+1}})\}$. This ensures
that properties (1$)_{i+1}$ and (6$)_{i+1}$ hold.

\vspace{4mm}

 Now is the time to choose compact polyhedra $P_{i+1}^\infty=P_{i+1}^{i+1},
P_{i+1}^i, \ldots,P_{i+1}^{\nu(i)},\ldots ,P_{i+1}^1$ in
$I^{n_{i+1}}$. First note that there is an open neighborhood
$\widetilde{V}$ of $p_{n_{i+1}}(X)$ in $I^{n_{i+1}}$ such that
$\widetilde{V} \times Q_{n_{i+1}} \subset N(X, \frac{2}{i+1})$.
Choose a compact polyhedron $P_{i+1}^\infty \subset I^{n_{i+1}}$
so that
\begin{enumerate}
\item[(17)] \
$p_{n_{i+1}}(X) \subset \int_{I^{n_{i+1}}} P_{i+1}^\infty
\subset P_{i+1}^\infty \subset \widetilde{V}$, \ and \
 $P_{i+1}^\infty \subset p_{n_i}^{-1}(\int_{I^{n_i}}(P_i^\infty)).$
\end{enumerate}
This can be done because (3)$_i$ implies (13)$_i^\infty$, i.e., $p_{n_i}(X)=p_{n_i}(p_{n_{i+1}}(X))\subset \int_{I^{n_i}}(P_i^\infty)$,
so\\ $p_{n_{i+1}}(X)\subset  p_{n_i}^{-1}(\int_{I^{n_i}}(P_i^\infty))$.
Note that (17) implies properties (3$)_{i+1}$ and (4$)_{i+1}$ for
$P_{i+1}^\infty$.
To satisfy the first part of $(2)_{i+1}$, we name $P_{i+1}^k=P_{i+1}^\infty$ for all $k\geq i+1$.

\vspace{2mm}

Let us now choose $P_{i+1}^k$, for $k=i, i-1,\ldots, 1$,  which we do by a downward recursion.

If $k > \nu(i)$, then here is how we make our choice: find an open neighborhood $\widetilde{V}_k$ of
$p_{n_{i+1}}(X_k)$ in $I^{n_{i+1}}$ such that $\widetilde{V_k}
\times Q_{n_{i+1}} \subset N(X_k, \frac{2}{i+1})$. Choose a
compact polyhedron $P_{i+1}^k \subset I^{n_{i+1}}$ so that
\begin{enumerate}
\item[(18)] \
$p_{n_{i+1}}(X_k) \subset \int_{I^{n_{i+1}}} P_{i+1}^k \subset
P_{i+1}^k \subset \widetilde{V_k} $,\ and \\ $P_{i+1}^k
\subset p_{n_i}^{-1}(\int_{I^{n_i}}(P_i^k))\ \bigcap \
\int_{I^{n_{i+1}}}(P_{i+1}^{k+1})$.
\end{enumerate}

This can be done because (3)$_i$ implies (13)$_i^k$, i.e., $p_{n_i}(X_k)=p_{n_i}(p_{n_{i+1}}(X_k))\subset \int_{I^{n_i}}(P_i^k)$,
so $p_{n_{i+1}}(X_k)\subset  p_{n_i}^{-1}(\int_{I^{n_i}}(P_i^k))$. Also note that
$p_{n_{i+1}}(X_k)\subset \int_{I^{n_{i+1}}}(P_{i+1}^{k+1})$, because before we reach the construction of $P_{i+1}^k$, $P_{i+1}^{k+1}$ is already constructed so that (13)$_{i+1}^{k+1}$ is true, so $p_{n_{i+1}}(X_{k+1}) \subset \int_{I^{n_{i+1}}}(P_{i+1}^{k+1})$, and also recall that $X_k\subset X_{k+1}\subset X$.

Note that (18) implies properties
   (2$)_{i+1}$ (the second part), (3$)_{i+1}$ and (4$)_{i+1}$ for
$P_{i+1}^k$, when $\nu(i)<k\leq i$.

\vspace {2mm}

For $k = \nu(i)$, we require the above mentioned properties and, additionally, that $P_{i+1}^{\nu(i)} \times Q_{n_{i+1}} \subset U$, where $U$ is the neighborhood of $X_{\nu(i)}$ indicated in (16).

For $k< \nu (i)$, proceed with the construction of $P_{i+1}^k$ as in the case of $i\geq k>\nu(i)$.
Conclude that properties (2)$_{i+1}$--(4)$_{i+1}$ are now true for all $k\in \{1,2,\ldots ,i\}\cup\{\infty\}$ for which they apply.

\vspace{4mm}

Let $\tilde{f}:= f''|_{P_{i+1}^{\nu(i)} \times Q_{n_{i+1}}} \circ
i: P_{i+1}^{\nu(i)} \to \EW(P_i^{\nu(i)},\Z/p ,l_{\nu(i)})$,
where $i: P_{i+1}^{\nu(i)} \to P_{i+1}^{\nu(i)} \times
Q_{n_{i+1}}$ is the inclusion.

\vspace{3mm}

Choose $\da_{i+1}$ and  triangulations $\widetilde\tau_{i+1}^\infty$ and $\tau_{i+1}^\infty$ for
$P_{i+1}^{\infty}$, which are also triangulating all $P_{i+1}^k$ for $k<i$
(where $\widetilde\tau_{i+1}^k := \widetilde\tau_{i+1}^\infty |_{P_{i+1}^k}$ and  $\tau_{i+1}^k := \tau_{i+1}^\infty |_{P_{i+1}^k}$), so that
(7$)_{i+1}$, (8$)_{i+1}$ and (9$)_{i+1}$ hold. Here is how this is
done: begin by taking a triangulation $\widetilde\tau_{i+1}^\infty$ of
$P_{i+1}^\infty$, which also triangulates all $P_{i+1}^k$, such
that $\mesh \widetilde\tau_{i+1}^\infty< \frac{\e_{i+1}}{2}$. The open stars in $\widetilde\tau_{i+1}^k$ of the
vertices of $\widetilde\tau_{i+1}^k$ form a cover $\scu_{i+1}^k=\{\St (v, \widetilde\tau_{i+1}^k) \vert \ v \in (\widetilde\tau_{i+1}^k)^{(0)} \}$ for $P_{i+1}^k$, where $k\in \{1,2,\ldots ,i\}\cup\{\infty\}$.


Note that for $x\in X$, $p_{n_{i+1}}(x)$ has to belong to some $\St(v, \widetilde\tau_{i+1}^\infty)$.
Then
 for any $y \in \St(v, \widetilde\tau_{i+1}^\infty)$, $\rho (y,p_{n_{i+1}}(x))\leq
2\mesh \widetilde\tau_{i+1}^\infty<\e_{i+1}$, so $\overline{\St}(v, \widetilde\tau_{i+1}^\infty) \subset
\overline{N}(p_{n_{i+1}}(x), \e_{i+1})\cap P_{i+1}^\infty$.  
Analogously, since for $x\in X_k$, $p_{n_{i+1}}(x)$ has to belong to some $\St(v, \widetilde\tau_{i+1}^k)\subset\St(v, \widetilde\tau_{i+1}^\infty) $, we get
$\overline{\St}(v, \widetilde\tau_{i+1}^k) \subset
\overline{N}(p_{n_{i+1}}(x), \e_{i+1})\cap P_{i+1}^k$,  for $k\in \{1,2,\ldots ,i\}$.

On the other hand, since $P_{i+1}^k$ is compact for $k\in \{1,2,\ldots ,i\}\cup\{\infty\}$, each cover  $\scu_{i+1}^k$ of $P_{i+1}^k$ has a Lebesgue number $\la_{i+1}^k$,  $k\in \{1,2,\ldots ,i\}\cup\{\infty\}$.
Thus it
is enough to pick a $\da_{i+1}>0$ such that
$$4\da_{i+1}\ < \min \left( \left\{ \la_{i+1}^k\ : \ k\in \{1,2,\ldots ,i\}\cup\{\infty\} \right\} \cup \left\{ \frac{4}{2^{n_{i+1}-1}}\right\}\right).$$
Now (7)$_{i+1}$ is satisfied. Also, for
any $x \in X_k$,
$\overline{N}(p_{n_{i+1}}(x),2\da_{i+1})\cap P_{i+1}^k$ is contained in some $\St(v, \widetilde\tau_{i+1}^k)$, for a vertex $v\in (\widetilde\tau_{i+1}^k)^{(0)}$. Pick one such star, and call its
closure $P_{x,i+1}^k$. Notice that $P_{x,i+1}^k$ is contractible.
Thus we get $(9)_{i+1}$ for $k<i+1$:
$$\overline{N}(p_{n_{i+1}}(x),2\da_{i+1})\cap P_{i+1}^k \subset P_{x,i+1}^k \subset
\overline{N}(p_{n_{i+1}}(x), \e_{i+1})\cap P_{i+1}^k.$$
Analogously, we get $(9)_{i+1}$ for $k=\infty$ and $x\in X$.
Finally,
choose a triangulation $\tau_{i+1}^\infty$ so that it refines
$\widetilde\tau_{i+1}^\infty$, and so that $(8)_{i+1}$ is satisfied.

\vspace{3mm}

Now that we have a triangulation for $P_{i+1}^\infty$, and
therefore for $P_{i+1}^{\nu(i)}$ too, take a cellular
approximation $\bar f_i :P_{i+1}^{\nu(i)} \to
\EW(P_i^{\nu(i)},\Z/p, l_{\nu(i)})$ of $\tilde f
:P_{i+1}^{\nu(i)} \to \EW(P_i^{\nu(i)},\Z/p, l_{\nu(i)})$.
Since $P_{i+1}^{\nu(i)}\times Q_{n_{i+1}} \subset U$, (16) is
valid for any $ u \in P_{i+1}^{\nu(i)}$, that is, $\omega \circ
f''(u,0)$ and $p_{n_i}(u,0)=p_{n_i}(u)$ belong to the $\tilde
\e$-neighborhood of the same principal simplex $\sa \in \tau_i^\infty$. We also know
that $\omega \circ f''(u,0)$ belongs to a principal simplex $\gamma$
which is a neighbor of $\sa$ (the choice of $\tilde \e$ makes sure
that $\gamma$ and $\sa$ are neighbors). Note that $\omega \circ
f''(u,0)= \omega \circ f''\circ i(u)=\omega\circ \tilde{f}(u) \in \gamma$. Now
$\omega\circ\bar{f}_i(u)$ also belongs to $\gamma$, because $\bar{f}_i$
is a cellular approximation of $\tilde{f}$, and the properties of the
Edwards--Walsh resolution $\omega$ guarantee that $\tilde{f}(u) \in
\omega^{-1}(\gamma)$ implies that $\bar{f}_i(u) \in
\omega^{-1}(\gamma)$. So we have found a simplex $\gamma$ of
$\tau_{i}^\infty$ such that $\omega\circ \bar{f}_i(u) \in \gamma$,
and $p_{n_i}(u)$ belongs to the $\tilde \e$-neighborhood of the
closed star of a vertex $v$ that is a common vertex of $\gamma$
and $\sa$. Therefore $\omega \circ\bar{f}_i: P_{i+1}^{\nu(i)}\ra P_{i}^{\nu(i)}$ and
$p_{n_i}|_{P_{i+1}^{\nu(i)}}: P_{i+1}^{\nu(i)}\ra P_{i}^{\nu(i)}$ are $\mathcal{W}$-near, and therefore  $\mathcal{V}$-homotopic. According to Lemma~\ref{L2.1} there exists
a continuous extension $\varphi :P_{i+1}^\infty \to P_{i}^\infty$ of $\omega\circ
\bar{f}_i$ such that $\vp$ and $p_{n_i}|_{P_{i+1}^\infty}$ are
$\mathcal{V}$-homotopic, and therefore $\mathcal{V}$-near.

\small{
\begin{displaymath}
\xymatrix{
&&& \EW(P_i^{\nu(i)},\Z/p, l_{\nu(i)})  \ar[d]^{\omega} &&\\
&&& P_{i}^{\nu(i)} \ar@{^{(}->}[rr]&&  \quad P_{i}^\infty\\
&&&&&\\
P_{x,i+1}^{\nu(i)} \ar@{^{(}->}[r] & P_{i+1}^{\nu(i)} \ar@{_{(}->}[d]^{i} \ar@/^/@{-->}[rruuu]^{\bar f_i} \ar[rruuu]|{\tilde f} \ar[rruu]|{\quad \qquad \quad p_{n_i}\vert \quad}  \ \  \ar@{^{(}->}[rr]&& P_{i+1}^\infty \ar[rruu]|{\quad p_{n_i}\vert \quad} \ar@/_/@{.>}[rruu]_{g_i^{i+1}}  \ar@/^/@{.>}[rruu]^{\vp}&&\\
U & P_{i+1}^{\nu(i)}\times Q_{n_{i+1}} \ar@{.>}[rruuuu]_{\!\!\!\!\!\! f''\vert} \ar@{_{(}->}[l] & X_{\nu(i)} \ar@{.>}[ruuu]|{p_{n_i}\vert } \ar@{_{(}->}[l] &&&
}
\end{displaymath}
}

\normalsize

With this, $(4)_{i+1}$, and the fact that we could have chosen
$\mathcal{V}$ as fine as we wish, we may assume that
$\vp(P_{i+1}^k)\subset P_i^k$, for all $1\leq k\leq\infty$.

\vspace{1mm}

Finally, making $\tau_{i+1}^\infty$ finer if necessary (but so that the properties of collapsibility required in (8)$_{i+1}$ are still preserved), take
$g_i^{i+1}:P_{i+1}^\infty \to P_{i}^\infty$  to be a simplicial
approximation of $\vp$. Therefore, for any $u \in P_{i+1}^\infty$,
there exists a simplex $\sa \in \tau_{i}^\infty$
 such that $g_{i}^{i+1}(u)$, $\vp(u) \in  \sa$. We also know that
 $\rho (\vp(u),p_{n_i}(u)) < \mesh \mathcal{V} < \frac{\da_i}{2}$,
 so $p_{n_i}(u) \in N(\sa ,\frac{\da_i}{2})$, i.e., property
 ($10)_{i+1}$ is true. Property ($11)_{i+1}$ is true because $g_i^{i+1}$ is a
simplicial approximation of $\vp$.

For property ($12)_{i+1}$, first notice that $g_i^{i+1}\vert_{P_{i+1}^{\nu(i)}} \simeq \vp\vert_{P_{i+1}^{\nu(i)}} = \omega \circ
\bar{f}_{i}$. Also, $\omega \circ \bar{f}_i$ and
$p_{n_i}|_{P_{i+1}^{\nu(i)}}$ being $\mathcal{W}$-near implies that
for all $x \in X_{\nu(i)}$, $\omega \circ \bar{f}_{i}
(P_{x,i+1}^{\nu(i)} ) \subset P_{x,i}^{\nu(i)}$. To see why, take any $u \in B_{x,i+1}^{\nu(i)\#}:= \overline
N(p_{n_{i+1}}(x),\e_{i+1})\cap P_{i+1}^{\nu(i)}$, i.e.,
$\rho (u,p_{n_{i+1}}(x)) < \e_{i+1}$;  by (5)$_{i+1}$,
$\rho(p_{n_i}(u),p_{n_i}(x))<\da_i$. Therefore, since $\mesh(\mathcal{W})<\frac{\da_i}{2}$,
 $$\rho(\omega \circ
\bar{f}_i(u),p_{n_i}(x))\leq \rho(\omega \circ
\bar{f}_i(u),p_{n_i}(u)) + \rho(p_{n_i}(u),p_{n_i}(x))<
\frac{\da_i}{2} + \da_i < 2\da_i,$$
so $\omega \circ \bar{f}_i(u)
\in B_{x,i}^{\nu(i)}:= \overline N(p_{n_{i}}(x),2\da_i)\cap P_i^{\nu(i)}$. Thus $\omega \circ \bar{f}_i
(B_{x,i+1}^{\nu(i)\#}) \subset B_{x,i}^{\nu(i)}$. Since $P_{x,i+1}^{\nu(i)}\subset\overline N(p_{n_{i+1}}(x),\e_{i+1})$, $\omega \circ
\bar{f}_{i} (P_{x,i+1}^{\nu(i)} ) \subset P_{x,i}^{\nu(i)}$, too.

Also, $\vp (P_{x,i+1}^{\nu(i)} )= \omega \circ \bar{f}_{i}
(P_{x,i+1}^{\nu(i)} ) \subset P_{x,i}^{\nu(i)}$, so $g_i^{i+1}$,
being a simplicial approximation of $\vp$, has the property
$g_i^{i+1}(P_{x,i+1}^{\nu(i)} ) \subset P_{x,i}^{\nu(i)}$.
Finally, $g_i^{i+1}\vert_{P_{x,i+1}^{\nu(i)}} \simeq \vp\vert_{P_{x,i+1}^{\nu(i)}}
=\omega \circ \bar{f}_{i}\vert_{P_{x,i+1}^{\nu(i)}}$, so property (12)$_{i+1}$ holds. $\square$


\begin{remark}
 Note that from our construction of $Z$, it follows that in general $Z$ is infinite dimensional.
\end{remark}

\begin{remark}
If we take $1<2< \ldots < m< \ldots$ instead of $l_1\leq l_2\leq \ldots \leq l_m\leq \ldots$, the Theorem~\ref{T1.1} becomes parallel to the result for $\dim_\Z$ from \cite{AJR}.\\
If $l_i=l_{i+1}$ but $X_{l_i}\subsetneq X_{l_{i+1}}$, we get $A_i=A_{i+1}$, but $Z_i\subsetneq Z_{i+1}$.
\end{remark}

What if the sequence of nonempty closed subspaces $X_1\subset X_2\subset\dots$
of the compact metrizable space $X$ from the statement of Theorem~\ref{T1.1} is finite, that is, we are given
$X_1\subset X_2\subset\dots \subset X_m\subset X$?
And what if $X$ itself is replaced by an $X_m$, i.e., we have $X_1\subset X_2\subset\dots \subset X_m=X$, where
for each $k\in \{1,2,\dots , m\}$, $\dim_{\Z/p} X_k\leq l_k$?

In either of these cases, Theorem~\ref{T1.1} yields a compact metrizable space $Z$ with closed subspaces
$Z_1\subset Z_2\subset\dots \subset Z_m\subset Z$, as well as a cell-like map $\pi: Z \ra X$ with all of the properties mentioned in Theorem~\ref{T1.1}, but we can adapt the proof so that it would use fewer polyhedra.

Namely, here are the changes that somewhat simplify the proof of Theorem~\ref{T1.1} in both of the finite cases mentioned above.

First, take a function $\nu :\N \ra \{1,2,\ldots ,m\}$ such that (i) and (ii) are still satisfied.

Second, change the conditions (2)$_{j\geq1}$ and (3)$_{j\geq1}$ from the original proof to the following:
\begin{enumerate}
\item[$(2)'_{j\ge 1}$] if $k\geq \min \ \{j, m+1 \}$ then $P_j^k=P_j^\infty$,
 and\\ $P_j^r\subset\int_{I^{n_j}}P_j^{r+1}$ whenever $r< \min \ \{j,m+1 \}$;
\item[$(3)'_{j\ge 1}$] $X\subset\int_Q(P_j^\infty\times Q_{n_j})
\subset N(X,\frac{2}{j})$, and,\\ whenever $k< \min \ \{j,m+1 \}$,
 $X_k\subset\int_Q(P_j^k\times Q_{n_j})\subset
N(X_k,\frac{2}{j})$;
\end{enumerate}

This will ensure that we produce only $m+1$ sequences of polyhedra $(P_j^k)_{j\in \N}$, $k \in \{1,\ldots ,m+1\}$, rather than the countably many sequences that were required in the original proof for $X_1\subset X_2\subset\dots \subset X_m\subset \dots \subset X$.

The rest of the proof is the same, provided that the change in indexes from $(2)'$  is taken into account in the remainder of the proof.

It is worth noting that, in the case when $X=X_m$,  the property $(3)'$ implies that we can take $P_j^\infty = P_j^m, \forall j$. Still, $Z$ and $Z_m$ would be different, since $Z_m=\lim((P_i^m)^{(l_m)},g_i^{i+1})$, and $Z=\lim(P_i^\infty,g_i^{i+1})=\lim(P_i^m,g_i^{i+1})$.
Also, the map $\pi\vert_{Z_m}: Z_m\ra X$ is a surjective $\Z/p$-acyclic map, while $\pi: Z \ra X$ is cell-like.

\begin{remark} \label{R-Dr}
In particular, for $m=1$ and $X=X_1$ such that $\dim_{\Z/p} X_1 \leq l_1$, Theorem~\ref{T1.1} produces a compact metrizable space $Z_1$ with $\dim Z_1 \leq l_1$, and a surjective $\Z/p$-acyclic map $\pi :Z_1 \ra X_1$.
So Theorem~\ref{T1.1} is indeed a generalization of Dranishnikov's resolution Theorem~\ref{Dr}.
\end{remark}

\section{Proof of a particular case of Theorem~1.1}

What follows is an outline of a proof for a particular case that Theorem~1.1 is covering, namely for the case when the sequence $l_1\leq l_2 \leq \ldots $ of upper bounds for $\dim_{\Z/p}$ does not become permanently stationary at any point. This proof was suggested to us by an anonymous referee. It does not work 
if this sequence is eventually constant, that is, if the spaces $X_i$ keep changing, but from some point $i_0$ on we have $l_{i_0}=l_{i_0+1}=\ldots$.

 For the sake of simplicity, let us suppose that $l_1<l_2<l_3\dots$ since the proof of this case can be adjusted to work for all cases in which the sequence is not eventually constant. 

Let $X_1\subset X_2 \subset \ldots$ be a sequence of nonempty closed subspaces of a compact metrizable space $X$ such that $\dim_{\Z/p} X_k \leq l_k$, $\forall k \in \N$. Apply Dranishnikov's Theorem~\ref{Dr} to $X_1$ in order to build a compact metrizable space $Z_1$ and a $\Z/p$-acyclic map $q_1:Z_1 \ra X_1$ such that $\dim Z_1\leq l_1$.
Let $Y_1=X\cup M(q_1)$ be the union of $X$ and the mapping cyllinder of $q_1$. Notice that the projection $p_1:Y_1\ra X$ is cell-like and that $\dim_{\Z/p} M(q_1) \leq l_1+1\leq l_2$, which makes $\dim_{\Z/p} X_2 \cup M(q_1) \leq l_2$. In order to produce $Z_2$ and $q_2$, apply Theorem~\ref{Dr} to  $X_2 \cup M(q_1)$, with the exception of requiring that $q_2$ has the property that $q_2|_{q_2^{-1}(Z_1)}$is a homeomorphism onto $Z_1$. 
Then put $Y_2=X\cup M(q_2)$ and $p_2: Y_2 \ra X$ to be the projection. Keep the procedure inductively and define $Z$ as the inverse limit of the inverse sequence $Y_1\leftarrow Y_2\leftarrow \ldots \leftarrow Y_k\leftarrow \ldots$


\end{document}